\documentclass[11pt]{article}

\newcommand{\documentdate}{December 12, 2017}

\usepackage{a4wide,epstopdf,latexsym,
            graphicx,varioref,color,amsmath,cite}
\DeclareGraphicsExtensions{.eps,.jpg}

\title{Optimality of orders one to three and beyond:\\
  characterization and evaluation complexity\\
  in constrained nonconvex optimization
  }
\author{C. Cartis\thanks{
   Mathematical Institute, 
   Oxford University,
   Oxford OX2 6GG, Great Britain.
   Email: coralia.cartis@maths.ox.ac.uk.}
~ N. I. M. Gould\thanks{
   Numerical Analysis Group,
   Rutherford Appleton Laboratory,
   Chilton OX11 0QX, Great Britain.
   Email:  nick.gould@stfc.ac.uk.}
~and Ph. L. Toint\thanks{
   Namur Center for Complex Systems (naXys) and Department of Mathematics,
   University of Namur,                
   61, rue de Bruxelles, B-5000 Namur, Belgium.
   Email: philippe.toint@unamur.be.
  }
}
\date{\documentdate}

\newcommand{\ass}[2]{\label{ass-#1}
                     \begin{list}{}{\setlength{\leftmargin}{2.5cm}}
                     \item \hspace{-2.5cm} \framebox[2.0cm]{\bf #1} \,\,#2
                     \end{list} }
\newcommand{\inneralg}{{Algorithm~\sc{inner}}}
\newcommand{\outeralg}{{Algorithm~\sc{outer}}}
\newcommand{\ep}{\epsilon_{\mbox{\tiny P}}}
\newcommand{\ed}{\epsilon_{\mbox{\tiny D}}}
\DeclareMathOperator*\globmin{globmin}

\newcommand{\numsection}[1]{\section{#1}\setcounter{equation}{0}}
\newcommand{\appnumsection}[1]{\section*{#1}\setcounter{equation}{0}
  \renewcommand{\theequation}{A.\arabic{equation}}
  \renewcommand{\thetheorem}{A.\arabic{theorem}}
  \renewcommand{\thetable}{A.\arabic{table}}
  \renewcommand{\thefigure}{A.\arabic{figure}} }
 
\newcommand{\calC}{{\cal C}} \newcommand{\calD}{{\cal D}}
 \newcommand{\calF}{{\cal F}}
 
\newcommand{\calI}{{\cal I}} 
\newcommand{\calK}{{\cal K}} 
\newcommand{\calM}{{\cal M}} \newcommand{\calN}{{\cal N}}
 \newcommand{\calP}{{\cal P}}
 
\newcommand{\calS}{{\cal S}} \newcommand{\calT}{{\cal T}}
 
 \newcommand{\calX}{{\cal X}}
 
\newcommand{\eqdef}{\stackrel{\rm def}{=}}
\newcommand{\tim}[1]{\;\; \mbox{#1} \;\;}
\newcommand{\beqn}[1]{\begin{equation}\label{#1}}
\newcommand{\eeqn}{\end{equation}}
\newcommand{\req}[1]{(\ref{#1})}
\newcommand{\ii}[1]{\{1, \ldots, #1 \}}
\newcommand{\iibe}[2]{\{ #1, \ldots, #2 \}}
\newcommand{\bctable}[1]{\begin{table}[htbp]
                         \begin{center}
                         \begin{tabular}{#1} }
\newcommand{\ectable}[1]{\end{tabular}
                         \caption{#1}
                         \end{center}
                         \end{table}}
\newlength{\thmw}
\newcommand{\lbthm}[3]{\vspace{\baselineskip}\noindent\hbox{%
  \lower\fboxrule\hbox{\vbox{\hrule\hbox{\vrule \kern-\fboxrule \vbox{%
  \vspace{\fboxsep} \noindent\hspace{2\fboxsep}\parbox{\thmw}{
  \begin{theorem}\label{#1}{\rm #2}\end{theorem}\vspace{-\lastskip}}
  \hspace{\fboxsep}}\kern-\fboxrule \vrule }}}}\newpage \hbox{%
  \lower\fboxrule\hbox{\vbox{\hbox{\vrule \kern-\fboxrule \vbox{%
  \noindent\hspace{2\fboxsep}\parbox{\thmw}{\rm #3}\hspace{\fboxsep}
  \vspace{4\fboxsep}}\kern-\fboxrule \vrule }\hrule }}}\vspace{\baselineskip}
}
\newtheorem{theorem}{Theorem}[section]
\newtheorem{lemma}[theorem]{Lemma}
\newtheorem{corollary}[theorem]{Corollary}

\newcommand{\ms}{\;\;\;\;}
\newcommand{\bpr}{{\bf Proof.} \hspace{1.5mm}}
\newcommand{\epr}{\hfill $\Box$ \vspace*{1em}}

\newcommand{\proof}[1]{
\begin{list}{}{
\setlength{\topsep}{0.0pt}
\setlength{\partopsep}{0.0pt}
\setlength{\leftmargin}{0.025\textwidth}
\setlength{\rightmargin}{0.5\leftmargin}
\setlength{\labelwidth}{0.5\leftmargin}
\setlength{\labelsep}{0.25\leftmargin}}
\item \bpr #1 \epr \noindent
\end{list}}
\newcommand{\half}{\sfrac{1}{2}}
\newcommand{\sfrac}[2]{{\scriptstyle \frac{#1}{#2}}}
\newcommand{\lcor}[2]{\vspace{\baselineskip} 
\noindent\framebox[\textwidth]{\parbox{0.95\textwidth}{
\begin{corollary} \label{#1} \rm #2 \end{corollary} } } \vspace{\baselineskip} }
\newcommand{\sixth}{\sfrac{1}{6}}
\newcommand{\smallRe}{\hbox{\footnotesize I\hskip -2pt R}}
\newcommand{\rvect}[1]{\left( \begin{array}{r} #1 \end{array} \right) }
\newcommand{\cvect}[1]{\left( \begin{array}{c} #1 \end{array} \right) }
\newcommand{\mat}[2]{\left(\begin{array}{#1}#2\end{array}\right)}
\newcommand{\spanset}[1]{\mbox{span}\left\{ #1 \right\}}
\newcommand{\s}[1]{^{\mbox{\protect\tiny #1}}}
\newcommand{\mystack}[2]{_{\stackrel{\scriptstyle #1}{\scriptstyle #2}}}
\newcommand{\kap}[1]{\kappa_{\mbox{\rm \tiny #1}}}
\newcommand{\bigmax}{\displaystyle \max}
\newcounter{algo}[section]
\renewcommand{\thealgo}{\thesection.\arabic{algo}}

\newcommand{\llem}[2]{\vspace{\baselineskip} 
\noindent\framebox[\textwidth]{\parbox{0.95\textwidth}{
\begin{lemma} \label{#1} \rm #2 \end{lemma} } } \vspace{\baselineskip} }
\newcommand{\lthm}[2]{\vspace{\baselineskip} 
\noindent\framebox[\textwidth]{\parbox{0.95\textwidth}{
\begin{theorem} \label{#1} \rm #2 \end{theorem} } } \vspace{\baselineskip} }
\newcommand{\Na}{\hbox{I\hskip -1.8pt N}}
\newcommand{\bigfrac}[2]{\frac{\displaystyle #1}{\displaystyle #2}}
\newcommand{\third}{\sfrac{1}{3}}
\newcommand{\bigsum}{\displaystyle \sum}

\newcommand{\balgo}[4]{\refstepcounter{algo}
\begin{center}\begin{figure}[htbp]
\vspace{\baselineskip}\noindent\hbox{%
  \lower\fboxrule\hbox{\vbox{\hrule\hbox{\vrule \kern-\fboxrule \vbox{%
  \vspace{\topsep} \noindent\hspace{2\fboxsep}\parbox{\thmw}{\vspace{0.5\topsep}
  {\bf Algorithm \thealgo : #2}\label{#1}\\
  \vspace*{-\topsep} \mbox{ }\\
  {\rm #3}\vspace{-\lastskip}}
  \hspace{\fboxsep}}\kern-\fboxrule \vrule }}}}
\end{figure}\end{center}
\newpage \hbox{%
  \lower\fboxrule\hbox{\vbox{\hbox{\vrule \kern-\fboxrule \vbox{%
  \noindent\hspace{2\fboxsep}\parbox{\thmw}{\rm #4}\hspace{\fboxsep}
  \vspace{\fboxsep}}\kern-\fboxrule \vrule }\hrule }}}\vspace{\baselineskip}
}

\begin{document}

\maketitle

\begin{abstract}
Necessary conditions for high-order optimality in smooth nonlinear constrained
optimization are explored and their inherent intricacy discussed. A two-phase
minimization algorithm is proposed which can achieve approximate first-,
second- and third-order criticality and its evaluation complexity is analyzed
as a function of the choice (among existing methods) of an inner algorithm for
solving subproblems in each of the two phases.  The relation between
high-order criticality and penalization techniques is finally considered,
showing that standard algorithmic approaches will fail if approximate
constrained high-order critical points are sought.
\end{abstract}

{\small
\textbf{Keywords:} nonlinear optimization, constrained problems, high-order optimality
conditions, complexity theory.
}
\vspace*{1cm}

\numsection{Introduction}

Analyzing the evaluation complexity of algorithms for solving the nonlinear
nonconvex optimization problem has been an active research area over the past
few years: we refer the interested reader to
\cite{AnanGe16,BergDiouGrat17,BianChen13,BianChen15,BianChenYe15,
  BirgGardMartSantToin17,BirgGardMartSantToin16,BoumAbsiCart16,
  CartGoulToin11,CartGoulToin11b,CartGoulToin11d,CartGoulToin12e,
  CartGoulToin13a,CartGoulToin10a,CartGoulToin12b,CartGoulToin15b,CartGoulToin16a,
  CartSampToin15,CartSche15,ChenToinWang17,CurtRobiSama17,CurtRobiSama17b,DodaViceZhan15,Duss15,
  GarmJudiVice15,GhadLan16,GeJianYe11,GrapYuanYuan15a,GrapYuanYuan15b,
  GratRoyeViceZhan15,GratRoyeVice17,GratSartToin08,Jarr13,LuWeiLi12,Mart17,MartRayd16,
  NestPoly06,NestGrap16,ScheTang14,UedaYama10,UedaYama10b,Vice13}
for contributions in this specific area.
The main focus of this thriving domain is to give (sometimes sharp) bounds on the
number of evaluations of a minimization problem's functions (objective and
constraints, if relevant) and their derivatives that are, in the worst case,
necessary for the considered algorithms to find an approximate critical point of a
certain order.  It is not uncommon that such algorithms involve costly
internal computations, provided the number of calls to the 
problem functions is kept as low as possible.

In nearly all cases, complexity bounds are given for the task of finding
$\epsilon$-approximate first- or (more rarely) second-order critical points,
typically using first- or second-order Taylor models of the objective function
in a suitable globalization framework such as those that use rust regions or
regularization. Notable exceptions are \cite{AnanGe16} where
$\epsilon$-approximate third-order critical points of unconstrained problems
are sought,
\cite{BirgGardMartSantToin17,CartGoulToin15a,CartGoulToin15b,
      CartGoulToin15c,BirgGardMartSantToin16} 
where $\epsilon$-approximate first-order critical points are considered using
Taylor models of order higher than two for unconstrained, convexly-constrained,
least-squares and equality-constrained problems, respectively, and
\cite{CartGoulToin16a} where general $\epsilon$-approximate $q$-th order
($q\geq 1$) critical points of convexly constrained optimization are analyzed
using Taylor models of degree $q$.

Because the present contribution focuses on problems involving a mixture of
convex inequality and nonlinear equality constraints, it is useful to set the
stage by considering earlier research in this constrained framework.  In
\cite{CartGoulToin12b}, the worst-case evaluation complexity of finding an
$\epsilon$-approximate first-order critical point for smooth nonlinear
(possibly nonconvex) optimization problems under convex constraints was
examined, using methods involving a second-order Taylor model of the objective
function. It was then shown that at most $O(\epsilon^{-3/2})$ evaluations of
the objective function and its derivatives are needed to compute such an
approximate first-order critical point.  This result, identical in order to
the best known result for the unconstrained case, assumes that the cost of
computing a projection onto the convex feasible set is neglible.  It comes
however at the price of potentially restrictive technical assumptions (see
\cite{CartGoulToin12b} for details). The analysis of \cite{CartGoulToin13a}
then built on this result by first specializing it to convexly constrained
nonlinear least-squares and then using the resulting complexity bound in the
context of a two-phase algorithm for a problem class involving general
constraints.  If $\ep$ and $\ed$ are the primal and dual criticality
thresholds, respectively, it was shown that at most $O(\ep^{-1/2}\ed^{-3/2})$
evaluations of the problem's functions and their gradients are needed to
compute an approximate critical point in that case, where the
Karush-Kuhn-Tucker (KKT) conditions are scaled to take the size of the
Lagrange multipliers into account. Because of the proof of this result is
based an the bound for the convex case, it suffers from the same limitations
(not to mention an additional constraint on the relative sizes of $\ep$ and
$\ed$, see \cite{CartGoulToin13a}). Another more general approach was
presented in \cite{Mart17} leading to the same complexity bounds, but at the
price of a subproblem involving the Jacobian of original nonlinear
constraints. The bounds derived in \cite{CurtRobiSama17b} for a trust-funnel
algorithm also consider a scaled KKT condition and are of the same order.  The
worst-case evaluation complexity of constrained optimization problems was also
recently analyzed in \cite{BirgGardMartSantToin16}, allowing for high-order
derivatives and models in a framework inspired by that of both
\cite{BirgGardMartSantToin17} and \cite{CartGoulToin12e,CartGoulToin13a}.  At
variance with these latter references, this analysis considers unscaled
approximate first-order critical points in the sense that such points satisfy
the standard unscaled KKT conditions with accuracy $\ep$ and $\ed$.  None of
these papers considers $\epsilon$-approximate second-order points for equality
constrained problems, except \cite{BoumAbsiCart16} where first- and
second-order optimality was proved for trust-region method defined on
manifolds.

The goal of this paper is twofold.  The first objective is to fill this gap by
deriving complexity bounds for finding $\epsilon$-approximate second- and third-order
critical points for the inequality/equality-constrained case.
The second is to examine higher-order optimality conditions (in the light of
\cite{CartGoulToin16a}) and to expose the intrinsic difficulties that arise
for criticality orders beyond three.

Our presentation is organized as follows. Necessary conditions for higher-order
criticality for nonlinear optimization problems involving both convex set
constraints and (possibly) nonlinear equality constraints are proposed and
discussed in Section~\ref{cnc-s}. A new two-phase algorithm is then introduced
in Section~\ref{general-s}, whose purpose is to compute $\epsilon$-approximate
critical points of orders one and two for such problems, and its evaluation
complexity is analyzed in Section~\ref{complexity-s} as a function of that of
an underlying inner algorithm for solving subproblems occuring in each of the
two phases. A discussion of the results and some conclusions are finally
presented in Section~\ref{disc-s}.

\noindent
\textbf{Basic notation.}
The notation in what follows is mostly inherited from
\cite{CartGoulToin16a}. $y^Tx$ denotes the Euclidean inner product of the
vectors $x$ and $y$ of $\Re^n$ and $\|x\| = (x^Tx)^{1/2}$ is the associated
Euclidean norm.
The cardinal of the set $\calS$ is denoted by $|\calS|$.
If $T_1$ and $T_2$ are tensors, $T_1\otimes T_2$ is their
tensor product and $\|T\|_q$ is the recursively induced Euclidean (or
spectral) norm of the $q$-th order tensor $T$.
If $T$ is a symmetric tensor of order $q$, the $q$-kernel of the
multilinear $q$-form
\[
T[v]^q \eqdef T[\underbrace{v, \ldots,v}_{q \tim{times}}]
\]
is denoted
$
\ker^q[T] \eqdef \{v \in \Re^n \mid T[v]^q = 0\}
$
(see
\cite{BrezTret03b,BrezTret05}).
Note that, in general, $\ker^q[T]$ is a union of cones\footnote{The
  1-kernels are not only unions of cones but also subspaces.
  However this is not true for general $q$-kernels, since both $(0,1)^T$ and
  $(1,0)^T$ belong to the 2-kernel of the non-negative symmetric 2-form
  $x_1x_2$ on $\Re^2$, but their sum does not. $\ker^1[x]$ is the usual
  orthogonal complement to the vector $x$, $\ker^2[M]$ is the standard
  nullspace of the matrix $M$.}.
If $\calX$ is a closed set, $\calX^0$ denotes its interior. The vectors
$\{e_i\}_{i=1}^n$ are the coordinate vectors in $\Re^n$. If $\{a_k\}$ and
$\{b_k\}$ are two infinite sequences of positive scalars converging to zero,
we say that $a_k= o(b_k)$ if and only if $\lim_{k \rightarrow \infty} a_k/b_k
= 0$. The normal cone to a general convex set $\calC$ at $x\in \calC$ is
defined by
\[
\calN_\calC(x) \eqdef \{ s \in \Re^n \mid s^T(z-x) \leq 0 \tim{for all} z \in
\calC \}
\]
and its polar, the tangent cone to $\calF$ at $x$,  by
\[
\calT_\calC(x)
= \calN_\calC^*(x)
\eqdef \{ s \in \Re^n \mid s^Tv  \leq 0 \tim{for all} v \in  \calN_\calC \}.
\]
Note that $\calC \subseteq \calT_\calC(x)$ for all $x \in \calC$.
We also define $P_{\calC}[\cdot]$ be the orthogonal projection onto $\calC$
and use the Moreau decomposition \cite{More62} which states that, for every
$x \in \calC$ and every $y \in \Re^n$ 
\beqn{moreau}
y = P_{\calT_\calC(x)}[y] + P_{\calN_\calC(x)}[y]
\tim{ and }
(P_{\calT_\calC(x)}[y]-x)^T(P_{\calN_\calC(x)}[y]-x) = 0.
\eeqn
(See \cite[Section 3.5]{ConnGoulToin00} for a brief
introduction of the relevant properties of convex sets and cones, or
\cite[Chapter 3]{HiriLema93} or \cite[Part I]{Rock70b} for an
in-depth treatment.)

\numsection{Necessary optimality conditions for constrained
  optimization}\label{cnc-s}

We  consider the  smooth constrained problem in the form
\beqn{genprob}
\min_{x \in \calF} f(x) 
\tim{ subject to }
c(x) = 0
\eeqn
where $c:\Re^n \to \Re^m$ is sufficiently smooth and $f$ and $\calF\subseteq \Re^n$ is
a non-empty, closed convex set. Note that this formulation covers the problems
involving both equality and inequality constraints, the latter being handled
using slack variables and the inclusion of the associated simple bounds in the
definition of $\calF$.

We start by investigating the necessary optimality conditions for problem
\req{genprob} at $x_*$ by considering possible feasible descent paths $x(\alpha)$ of
the form 
\beqn{feasarc}
x(\alpha) = x_* + \sum_{i=1}^q\alpha^i s_i + o(\alpha^q)
\eeqn
where $\alpha > 0$.  As in \cite{CartGoulToin16a}, we define the $q$-th order
descriptor set of $\calF$ at $x$ by 
\beqn{Dq-def}
\calD_\calF^q(x)
\eqdef \bigcup_{\varsigma >0}\left\{ (s_1, \ldots, s_q ) \in \Re^{n \times q} \mid x +
\sum_{i=1}^q\alpha^i s_i + o(\alpha^q)\in \calF  \right\}
\eeqn
Note that $\calD_\calF^1(x) = \calT_\calF(x)$, the standard tangent cone
to $\calF$ at $x$. We say that a feasible curve\footnote{Or arc, or path.}
$x(\alpha)$ is tangent to $\calD_\calF^q(x)$ if \req{feasarc} holds for some
$(s_1, \ldots, s_q) \in \calD_\calF^q(x)$.

The necessary optimality conditions for problem \req{genprob} also involve the
associated Lagrangian function 
\beqn{Lag-def}
\Lambda(x,y) \eqdef f(x) + y^Tc(x),
\eeqn
the subspace
\beqn{Mdef}
\calM(x)
\eqdef \ker^1[\nabla_x^1c(x)] \cap \ker^1[\nabla_x f(x)]
\eeqn
and the index sets
$\calP(j,k)$ defined, for $k \leq j$, by
\beqn{Pjks-def}
\calP(j,k) \eqdef \{ (\ell_1, \ldots ,\ell_k) \in \ii{j}^k
                    \mid \sum_{i=1}^k \ell_i = j \}.
\eeqn
For $k \leq j \leq 4$, these are given by Table~\ref{Pjks-t}.

\bctable{c|llll}
$j\downarrow$ & $k \rightarrow$ &&&\\
&
\multicolumn{1}{c}{1} &
\multicolumn{1}{c}{2} &
\multicolumn{1}{c}{3} &
\multicolumn{1}{c}{4} \\
\hline
1 & \{(1)\} & & & \\
2 & \{(2)\} & \{(1,1)\} & & \\
3 & \{(3)\} & \{(1,2),(2,1)\} & \{(1,1,1)\}  & \\
4 & \{(4)\} & \{(1,3),(2,2),(3,1))\} & \{(1,1,2),(1,2,1),(2,1,1)\} & \{(1,1,1,1)\}\\
\ectable{The sets $\calP(j,k)$ for $k\leq j \leq 4$ \label{Pjks-t}}

\lbthm{ccnc-th}{
Suppose that $f$ and each of the $\{c_i\}_{i=1}^m$ are $q$ times continuously
differentiable in an open set containing $\calF$, and that $x_*$ is a local 
minimizer for problem \req{genprob}.  
\noindent
Then we have that $c(x_*) = 0$ and, for some $y_* \in \Re^m$ and 
$j \in \iibe{1}{q}$,
\beqn{ccncq-e}
    \sum_{k=1}^j \frac{1}{k!}\left( \sum_{(\ell_1, \ldots,\ell_k) \in
     \calP(j,k)} \nabla_x^k \Lambda(x_*,y_*)[s_{\ell_1}, \ldots,
      s_{\ell_k}]\right)
    \geq 0
\eeqn
}{
for all $\{s_i\}_{i=1}^j$ such that $s_1 \in \calT_*$, $x(\alpha) \in
\calF$ for $\alpha >0$ sufficiently small, and such that
\beqn{ccncqz-e}
    \sum_{k=1}^i \frac{1}{k!}\left( \sum_{(\ell_1, \ldots,\ell_k) \in
    \calP(i,k)} \nabla_x^k \Lambda(x_*,y_*)[s_{\ell_1}, \ldots,
      s_{\ell_k}]\right)
    = 0,
\ms (i = 1,\ldots,j-1),
\eeqn
and
\beqn{ccncf}
\sum_{k=1}^i \frac{1}{k!} \left(\sum_{(\ell_1, \ldots,\ell_k) \in \calP(i,k)}
\nabla_x^k c(x_*)[s_{\ell_1}, \ldots, s_{\ell_k}]\right) = 0,
\ms (i = 1,\ldots,j).
\eeqn
}

\proof{
  Consider feasible paths of the form \req{feasarc}. Substituting this
  relation in the expression $f(x(\alpha)) \geq f(x_*)$ (which must be true
  for small $\alpha > 0$ if $x_*$ is a local minimizer) and collecting terms
  of equal degree in $\alpha$,  we obtain that, for sufficiently small $\alpha$,
  \beqn{fexpq}
  0
  \leq f(x(\alpha)) - f(x_*)
  = \sum_{j=1}^q \alpha^j
          \sum_{k=1}^j \frac{1}{k!} \Bigg(\sum_{(\ell_1, \ldots,\ell_k) \in \calP(j,k)}
          \nabla_x^k f(x_*)[s_{\ell_1}, \ldots, s_{\ell_k}]\Bigg)
    + o(\alpha^q)
  \eeqn
  where $\calP(i,k)$ is defined in \req{Pjks-def}.
  Similarly, substituting \req{feasarc} in the expression $c(x(\alpha)) = 0$ and
  collecting terms of equal degree in $\alpha$, we obtain that, for
  sufficiently small $\alpha$,
  \beqn{cexpq}
  0 = c(x(\alpha)) = \sum_{j=1}^q \alpha^j
          \sum_{k=1}^j \frac{1}{k!} \Bigg(\sum_{(\ell_1, \ldots,\ell_k) \in \calP(j,k)}
          \nabla_x^k c(x_*)[s_{\ell_1}, \ldots, s_{\ell_k}]\Bigg)
    + o(\alpha^q);
  \eeqn
     Adding now $f(x(\alpha))$ from \req{fexpq} to $y_*^Tc(x(\alpha))$ from
  \req{cexpq}, we obtain that
  \beqn{Lagcond}
  \sum_{j=1}^q \alpha^j \sum_{k=1}^j\frac{1}{k!}
  \Bigg( \sum_{(\ell_1, \ldots,\ell_k) \in \calP(j,k)} 
  \nabla_x^k \Lambda(x_*,y_*)[s_{\ell_1}, \ldots, s_{\ell_k}] \Bigg)
  + o(\alpha^q)
  \geq  0
  \eeqn
  for $\alpha >0 $ sufficiently small. For this to be true, we need each
  coefficient of $\alpha^j$ to be non-negative on the zero set of the
  coefficients $1, \ldots, j-1$ (i.e., satisfying \req{ccncqz-e}), subject to
  the requirement that the arc \req{feasarc} must be feasible for $\alpha$
  sufficiently small, that is \req{ccncf} holds and $x(\alpha))\in \calF$ for
  sufficiently small $\alpha > 0$.

  We start by examining first-order conditions ($q=1$). For $j=1$ (for which
  conditions \req{ccncqz-e} and \req{ccncf} are void) and observing that
  $\calP(1,1) = \{(1)\}$ (see Table~\ref{Pjks-t}), the necessary positivity of
  the coefficient of $\alpha$ in \req{Lagcond} implies that, for $s_1 \in \calT_*$,
  \beqn{nec1}
  \nabla_x^1\Lambda(x_*,y_*)[s_1]  \geq 0.
  \eeqn

  Consider now the case where $q=2$ and assume that $s_1 \in \calT_*$ and also that
  \req{ccncqz-e} and \req{ccncf} hold.  The former condition requires that $s_1
  \in \ker^1[\nabla_x^1c(x*)]$ and the latter that $s_1 \in
  \ker^1[\nabla_x^1\Lambda(x_*,y_*)]$, yielding together that 
  \[
  s_1 \in \calT_* \cap \ker^1[\nabla_x^1c(x_*)] \cap \ker^1[\Lambda(x_*,y_*)]
  =  \calT_* \cap \calM(x_*).
  \]
  Then the  coefficient of $\alpha^2$ in \req{Lagcond} must be non-negative,
  which yields,  using $\calP(2,1) =  \{(2)\}$, $\calP(2,2)= \{(1)\}$ (see
  Table~\ref{Pjks-t}) and \req{mgradLinN}, that  
  \beqn{ccnc2-e}
  \nabla_x^1\Lambda(x_*,y_*)[s_2] + \half \nabla_x^2\Lambda(x_*,y_*)[s_1]^2
  \geq 0.
  \eeqn
  which is \req{ccncq-e} for $q=2$.

  We may then proceed in the same manner for higher orders, each time
  considering them in the zero set of the previous coefficients  (that is
  \req{ccncqz-e}), and verify that \req{Lagcond} directly implies \req{ccncq-e}.
} 

\noindent
We note that, as the order $j$ grows,
\req{ccncq-e} and \req{ccncf} for $i=j$ may be interpreted as imposing
conditions on $s_j$ (via $\nabla_x^1\Lambda(x_*,y_*)[s_j]$ and
$\nabla_x^1f(x_*)[s_j]$), given the directions $\{s_i\}_{i=1}^{j-1}$ satisfying
\req{ccncqz-e} and \req{ccncf} for $i \in \ii{j-1}$.

Theorem \ref{ccnc-th} covers some well-known cases, as shown by the next
corollary.

\lcor{NCq1-c}{
Suppose that $f$ and each of the $\{c_i\}_{i=1}^m$ are $q$ times continuously
differentiable in an open set containing $\calF$, and that $x_*$ is a local 
minimizer for problem \req{genprob}. Let $\calN_*$ be the normal cone
to $\calF$ at $x_*$ and $\calT_*$ the corresponding tangent cone. 
\noindent
Then we have that $c(x_*) = 0$ and, for some $y_* \in \Re^m$,
\beqn{mgradLinN}
- \nabla_x^1\Lambda(x_*,y_*)
\in \calN_*,
\eeqn
Moreover,  if $x_* \in \calF^0$, the interior of $\calF$, then $\nabla_x^2
\Lambda(x_*,y_*)$ is positive semi-definite on $\ker^1[\nabla_x^1c(x_*)]$.
}

\proof{
  Using the fact that the normal cone $\calN_*$ is the polar of $\calT_*$, we
  immediately deduce from \req{nec1} that \req{mgradLinN} holds.
  If we also assume that $x_* \in \calF^0$, \req{mgradLinN} unsurprisingly
  reduces to $\nabla_x^1 \Lambda(x_*,y_*) = 0$, while, for $j=q=2$,
  \req{ccncq-e} gives that $\nabla_x^2 \Lambda(x_*,y_*)$ must be positive
  semi-definite on the subspace defined by \req{ccncf}, that is $\calM(x_*) =
  \ker^1[\nabla_x^2c(x_*)]$.
} 

\begin{figure}[htbp]
\centerline{
\includegraphics[height=6cm]{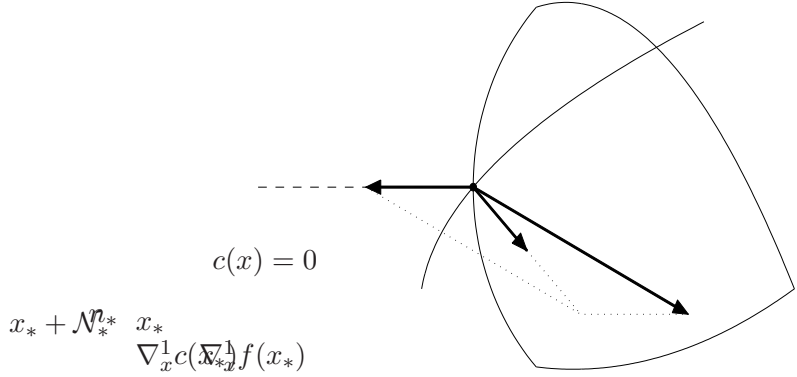}  
}
\begin{picture}(120,0.1)(0,0)
\put(25,28){$x_*+\calN_*$}
\put(55,32){$n_*$}
\put(73,29){$x_*$}
\put(97,16){$\nabla_x^1 f(x_*)$}
\put(73,16){$\nabla_x^1 c(x_*)$}
\put(95,40){$\calđ$}
²\put(102,52){$c(x)=0$}
\end{picture}
\caption{\label{bao-mx1rst_fig} The condition \req{mgradLinN} with $\calN_*$ shown
  as a dashed half line. Note that $n_* = -\nabla_x^1 \Lambda(x_*,y_*) \neq
  P_{\calN_*}[-\nabla_x^1f(x_*]$ (adapted from \cite{ConnGoulToin00}). }
\end{figure}

The conditions stated in Corallary~\ref{NCq1-c} for $q=1$ or $2$ are standard
  (for \req{mgradLinN}, see \cite[Theorem~3.2.1, p. 46]{ConnGoulToin00}, for
  instance, and Figure~\ref{bao-mx1rst_fig} for an illustration).   For more
  general cases, the complicated conditions \req{ccncq-e}-\req{ccncf} appear not
  to have been stated before and merit some discussion.

It was observed in \cite[Section~3]{CartGoulToin16a} that the necessary
optimality condition for the essentially unconstrained case where $x_* \in
\calF^0$ (implying $\calN_* = \{0\}$) combines more than a single derivative tensor
and $s_i$ for orders four and above.  If equality constraints are present this
situation already appears at order three (and above). Indeed, it can be verified
that the necessary conditions \req{ccncq-e} and \req{ccncf} for $q=3$ and
$\calN_* = \{0\}$ (and hence $\nabla_x^1\Lambda(x_*,y_*)=0$ because of
\req{mgradLinN}) can be written as 
\beqn{ccnc3}
\nabla_x^2\Lambda(x_*,y_*)[s_1,s_2] + \sixth \nabla_x^3\Lambda(x_*,y_*)[s_1]^3
= 0
\eeqn
for all $s_1 \in \calT_* \cap \ker^1[\nabla_x^1\Lambda(x_*,y_*)]
\cap \ker^2[\nabla_x^2\Lambda(x_*,y*)]$ and
\beqn{ccncf3}
\nabla_x^1c(x_*)[s_2] + \half \nabla_x^2c(x_*)[s_1]^2 = 0,
\ms
\nabla_x^1c(x_*)[s_3] + \nabla_x^2c(x_*)[s_1,s_2] + \sixth \nabla_x^3c(x_*)[s_1]^3= 0.
\eeqn
These conditions do not require the second term of the
left-hand side  of \req{ccnc3} to vanish.  This is at variance with the unconstrained case, since
  second-order necessary conditions then ensure that
  $\nabla_x^2\Lambda(x_*,y_*)$ is positive semidefinite on $\Re^n$ and therefore
  admits a square root. Thus $\nabla_x^2\Lambda(x_*,y_*)[s_1,s_2]=
  [\nabla_x^2\Lambda(x_*,y_*)^{\half}s_2]^T[\nabla_x^2\Lambda(x_*,y_*)^{\half}s_1]
  = 0$ since $s_1$ must belong to
  $\ker^2[\nabla_x^2\Lambda(x_*,y_*)]$. However, this
  argument no longer applies in the constrained case because
  $\nabla_x^2\Lambda(x_*,y_*)$ is only positive semidefinite on a strict subspace
  of $\Re^n$ and the square root may fail to exist, as is illustrated by the following
example.

\noindent
{\bf Example.} Consider the problem
\[
\min_{x \in \smallRe^3} x_1 + x_2^2 + x_2^3 -x_3
\tim{ subject to }
c(x) = \rvect{- x_1 - x_2^2 + x_1x_2 + x_3 \\
                x_1 + x_2^2 + x_1x_2 + x_3 } = 0,
\]
for which the origin is a high-order saddle point.

\begin{figure}[htbp]
\begin{center}
\vspace*{1.5mm}
\includegraphics[height=4.3cm]{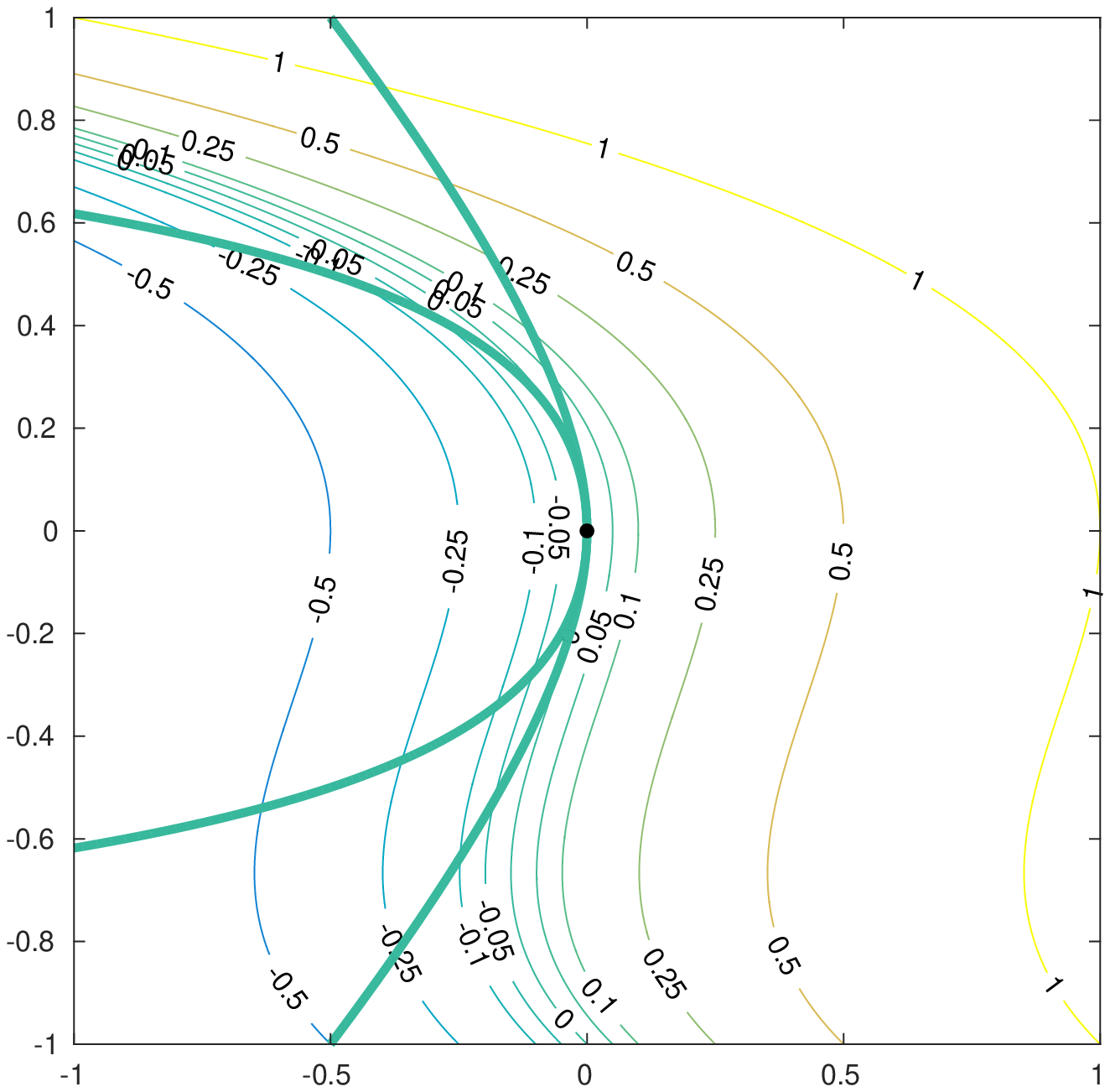}  
\hspace*{10mm}
\includegraphics[height=4.3cm]{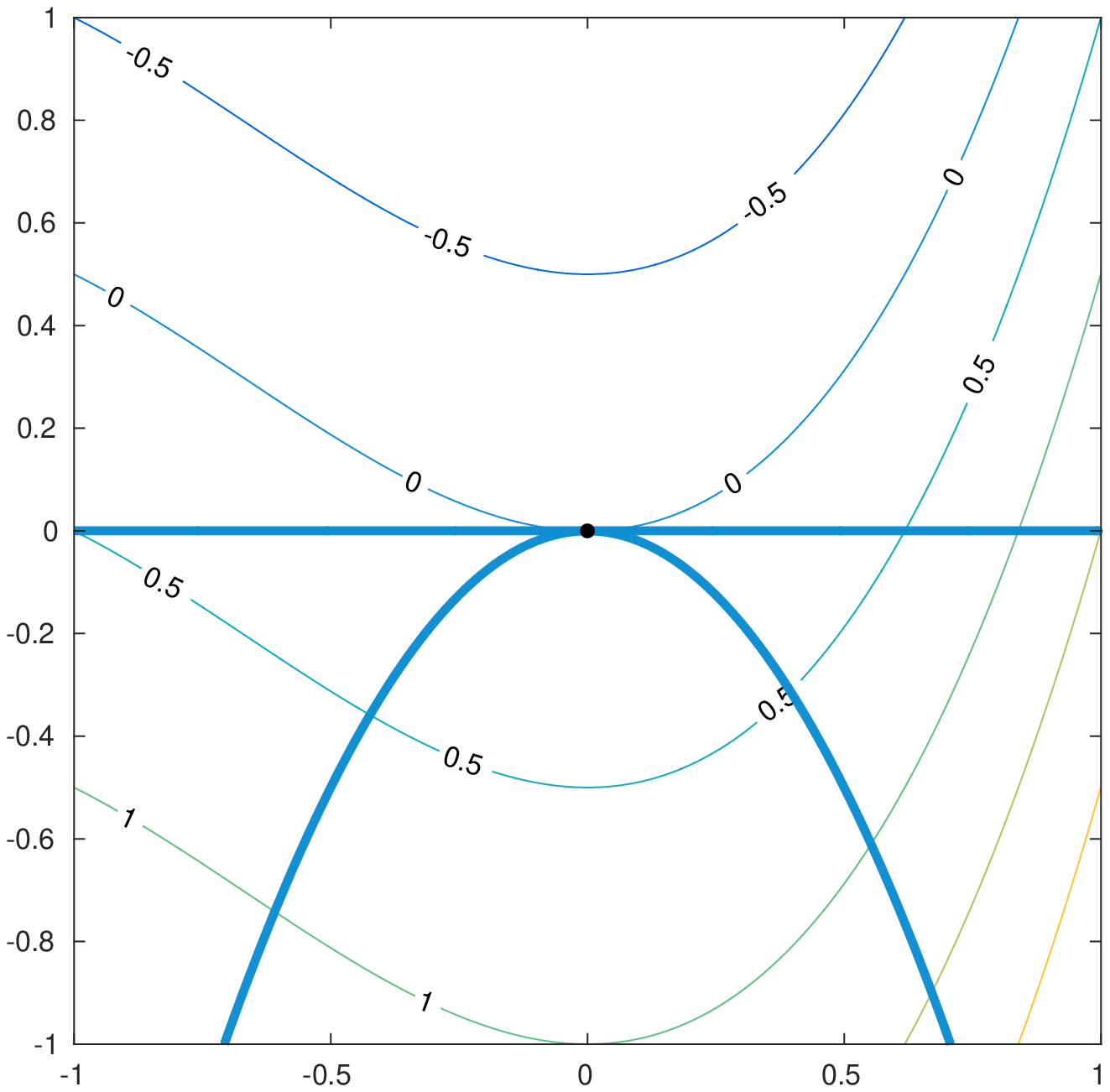}  
\hspace*{10mm}
\includegraphics[height=4.3cm]{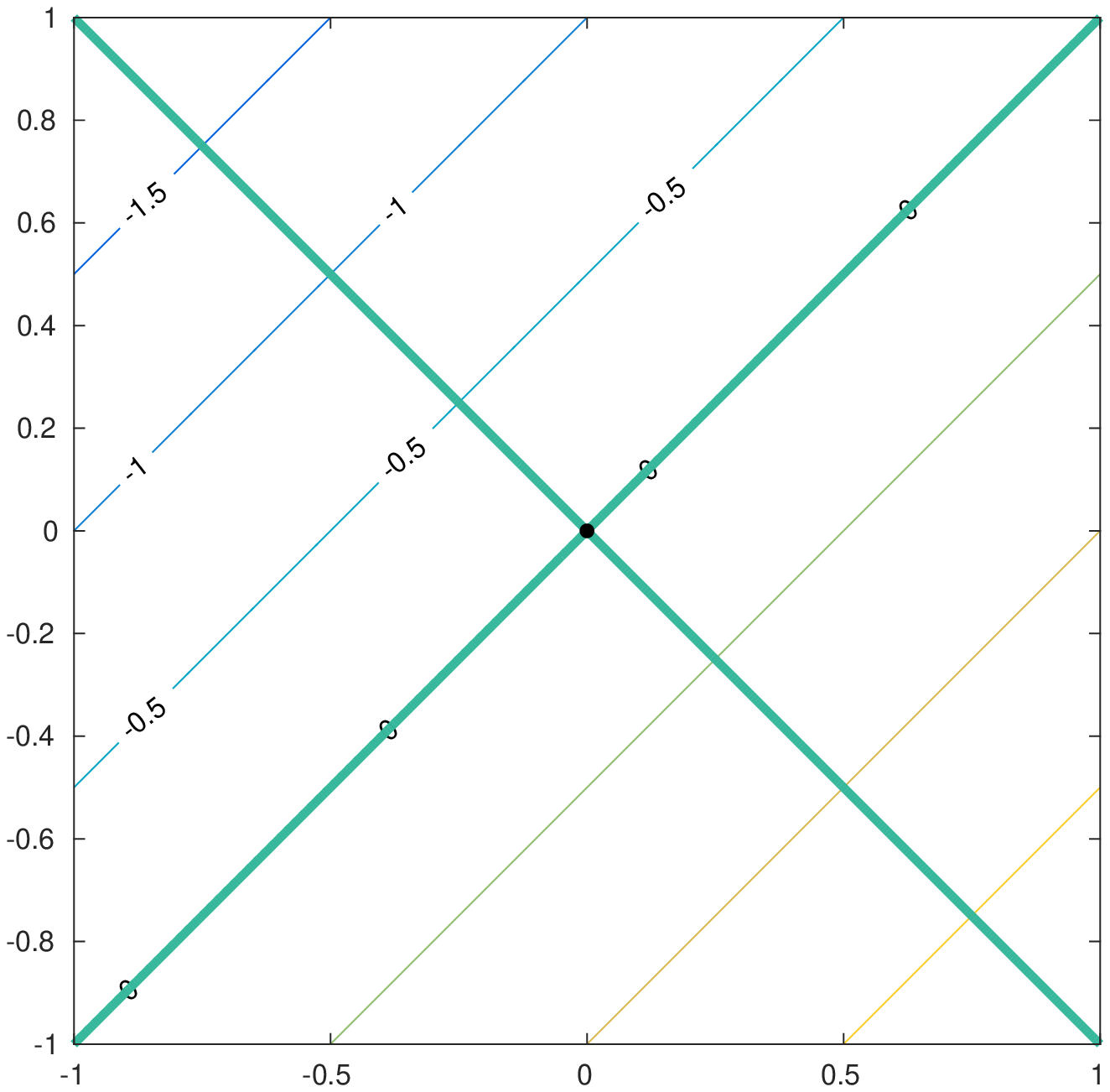}  
\caption{\label{ex4_fc}The contour lines of $f(x_1,x_2,0)$ (left)
  $f(0,x_2,x_3)$ (center), $f(x_1,0,x_3)$ (right) and the
  two constraints intersecting at the origin (thick).}
\end{center}
\end{figure}

Comparing the constraints' expression with \req{feasarc} for $q=3$, we see
that \req{Dq-def} holds for 
\[
s_1 = e_2,
\ms
s_2 = - e_1
\tim{ and }
s_3 = e_3
\]
since then
\[
x(\alpha) = \rvect{ -\alpha^2 \\ \alpha\;\, \\ \alpha^3 }
\tim{ and }
c(x(\alpha)) = \rvect{ \alpha^2 - \alpha^2 - \alpha^3 + \alpha^3 \\
                      -\alpha^2 + \alpha^2 - \alpha^3 + \alpha^3 }= 0.
\]
Now,
\[
\nabla_x^1f(x) = \cvect{1 \\ 2x_2 + 3 x_2^2 \\ -1 }
\ms
\nabla_x^2f(x) = \mat{ccc}{ 0 & 0      & 0 \\
                            0 & 2+6x_2 & 0 \\
                            0 & 0      & 0 }
\tim{ and }
[\nabla_x^3f(x)]_{2,2,2} = 6.
\]
\[
\nabla_x^1c(x) = \mat{rrr}{ -1 + x_2 & -2x_2 + x_1 & 1 \\
                             1 + x_2 &  2x_2 + x_1 & 1 },
\ms
\nabla_x^2c_1(x) = \mat{rrr}{ 0 & 1 & 0 \\
                             1 & -2 & 0 \\
                             0 &  0 & 0},
\ms
\nabla_x^3c_1(x) = 0,
\]
\[
\nabla_x^2c_2(x) = \mat{ccc}{ 0 & 1 & 0 \\
                             1 & 2 & 0 \\
                             0 & 0 & 0 }
\tim{ and }
\nabla_x^3c_2(x) = 0.
\]
Moreover,
\[
\begin{array}{ll}
\nabla_x^1 &\!\!\!\!c(0)[s_2] + \half \nabla_x^2c(0)[s_1]^2\\
  & =  -\mat{rrr}{ -1 & 0 & 1 \\ 1 &  0 & 1 }e_1
  + \half  \left[e_2^T\mat{rrr}{ 0 & 1 & 0 \\ 1 & -2 & 0 \\ 0 & 0 & 0 }e_2\right]e_1
  + \half  \left[e_2^T\mat{rrr}{ 0 & 1 & 0 \\ 1 &  2 & 0 \\ 0 & 0 & 0 }e_2\right]e_2\\
  & = 0.
\end{array}
\]
and
\[
\begin{array}{ll}
  \nabla_x^1&\!\!\!\!c(0)[s_3]+\nabla_x^2c(0)[s_1,s_2]
              +\sixth\nabla_x^3c(0)[s_1]^3\\*[2ex]
& = \mat{rrr}{ -1 & 0 & 1 \\ 1 &  0 & 1 }e_3
- \left[e_2^T\mat{rrr}{ 0 & 1 & 0 \\ 1 & -2 & 0 \\ 0 & 0 & 0 }e_1\right]e_1
- \left[e_2^T\mat{rrr}{ 0 & 1 & 0 \\ 1 &  2 & 0 \\ 0 & 0 & 0 }e_1\right]e_2
- \sixth 0^T [e_1]^3\\
& = 0.
\end{array}
\]
Thus \req{ccncf3} holds. From the values of $\nabla_x^1f(0)$ and
$\nabla_x^1c(0)$, we verify that setting $y_0 = (1,0)^T$ ensures
that $\nabla_x^1 \Lambda(0,y_0) = 0$. Hence \req{mgradLinN}  holds as well.
Moreover, we have that
\[
\ker^1[\nabla_x^1c(0)]
= \ker^1\left[\mat{rrr}{ -1  &  0 & 1 \\
                          1  &  0 & 1 }\right]
= \spanset{e_2},
\;\;
\nabla_x^2 \Lambda(0,y_0) = \left(\begin{array}{ccc}
  0 & 1 & 0 \\
  1 & 0 & 0 \\
  0 & 0 & 0
\end{array}\right)
\]
and the only nonzero component of $\nabla_x^3\Lambda(0,y_0)$ is its (2,2,2)
element which is 6.  Thus \req{ccncqz-e} also holds for $i=2$ .  In addition,
it is easy to check that the third-order necessary condition \req{ccnc3} holds with
\[
\nabla_x^2\Lambda(0,y_0)[s_1,s_2] = -1
\tim{ and }
\nabla_x^3\Lambda(0,y_0)[s_1]^3 = 6.
\]
This shows that the term involving $\nabla_x^3\Lambda(0,y_0)[s_1]^3$  is not the
only one occuring in the third-order necessary condition for our example
problem, as announced. Figure~\ref{ex4_fc} show the level lines of the objective
function and the constraint manifold in the $(x_1,x_2)$ $(x_2,x_3)$
$(x_1,x_3)$ planes, 
illustrating the interaction of the objective function's curvature and
feasible set. \epr

The third order necessary condition therefore must consider both terms in
\req{ccnc3} and cannot rely only on the third derivative of the Lagrangian
along a well-chosen direction or subspace.  In general, the $q$-th order
necessary conditions will involve (in \req{ccncq-e}) a mix of other terms than
those involving the $q$-th derivative tensor of the Lagrangian applied on
vectors $s_i$ for $i>1$, themselves depending on the geometry of the set of
feasible arcs. At this stage, for lack of a suitable formal understanding of
this geometry, conditions \req{ccncq-e}-\req{ccncf} remain very difficult to
interpret or check.

\numsection{A minimization algorithm}\label{general-s}

Having analyzed the necessary condition for problem \req{genprob} and seen
that conditions for orders above three are, at this stage, very difficult to
verify for general problems, we now describe a two-phase algorithm whose purpose is to
find approximate critical points of order one and two (and possibly three as
we discuss below). Since the presentation is independent of the order $q$ of the
critical points sought, we keep this order general in what follows.

\subsection{Inner algorithms for constrained least-squares problems}
\label{inner-s}

As was the case in \cite{CartGoulToin11b,CartGoulToin15a}, the new two-phase
algorithm relies on an inner algorithm for solving the convexly constrained
nonlinear least-squares problem in each of its phases.  We therefore start by
reviewing the existence and properties of algorithms for solving this
subproblem.

Consider first the standard convexly constrained problem
\beqn{ccprob}
\min_{x \in \calF} \psi(x)
\eeqn
where $\psi$ is a smooth function from $\Re^n$ to $\Re$ and $\calF$ is (as in
\req{genprob}) a non-empty closed convex set.  Following
\cite{CartGoulToin16a}, an $\epsilon$-approximate $q$-th order critical point
for this problem can be defined as a point $x$ such that
\beqn{cc-approx}
\phi_{\psi,j}^\Delta(x) \leq \ed\s{LS} \Delta^j
\tim{for}
j = 1,\ldots, q
\eeqn
and some $\Delta \in (0,1]$, where,
for $\calF(x) \eqdef \{ d \in \Re^n | \mid  x + d \in \calF \}$,
\beqn{phidef}
\phi_{\psi,j}^\Delta(x)
\eqdef \psi(x)-\globmin_{\mystack{d \in \calF(x)}{\|d\|\leq\Delta}}T_{\psi,j}(x,d),
\eeqn
is the largest feasible decrease of the $j$-th order Taylor model $T_{\psi,j}(x,s)$
achievable at distance at most $\Delta$ from $x$. Note that
$\phi_{\psi,j}^\Delta(x)$  is a continuous function of $x$ and $\Delta$ for given
$\calF$ and $f$ (see \cite[Theorem~7]{Hoga73}). It is also monotonically
increasing in $\Delta$. Also note that the global minimization involved in
\req{phidef} is efficiently solvable for $j=1$  because it is convex.  It is
also tractable in the unconstrained case for $j=2$ since it then reduces to a
trust-region subproblem.

Algorithms for finding $\epsilon$-approximate first-order critical points
for problem \req{ccprob}, i.e. points satisfying \req{cc-approx} for some
algorithm-dependent $\Delta \in (0,1]$ have already
been analyzed, for instance in \cite{CartGoulToin12b,CartGoulToin15b} or
\cite{CartGoulToin16a}, the first two being of the regularization type, the
last one being a trust-region method.  Such algorithms generate a sequence of
feasible iterates $\{x_k\}$ with monotonically decreasing objective-function
values $\{\psi(x_k)\}$. The method described in \cite{CartGoulToin15b} proceeds
by approximately minimizing models based on the regularized Taylor series of
degree $p$ and and it can be shown
\cite[Lemma~2.4]{CartGoulToin15b}\footnote{Observe that
  $\phi_{\psi,1}^\Delta(x)/\Delta =\chi_{\psi,1}(x)$ as defined in
  \cite[equation (2.4)]{CartGoulToin15b}, irrespective of the
  value of $\Delta \in (0,1]$.}
that, as long as the stopping criterion \req{cc-approx} fails for $q=1$ and
$\Delta=1$, a sufficient objective-function decrease
\beqn{cc-decr-succ-1}
\psi(x_k) - \psi(x_{k+1}) \geq \kap{decr}^\psi \epsilon^{\sfrac{p+1}{p}}
\eeqn
holds for each $k \in \calS$, where $\kap{decr}^\psi\in (0,1)$ is a
constant independent of $\epsilon$, and where $\calS$ is the set of
``successful iterations'' at which an effective step is made
(i.e.\ $x_{k+1}\neq x_k$).  Moreover, it can also be shown
\cite[Lemma~2.1]{CartGoulToin15b} that the set $\calS$ cannot be too small in
the sense that, for all $k\geq 0$,  
\beqn{cc-unsucc-bound}
k \leq \kap{uns}^\psi |\calS \cap \ii{k}|
\eeqn
for some constant $\kap{uns}^\psi>0$.  Both $\kap{decr}^\psi$ and
$\kap{uns}^\psi$ typically depend on the details of the considered algorithm
and of the Lipschitz constant associated with the highest derivative used in
the objective-function's model. Both \req{cc-decr-succ-1} and
\req{cc-unsucc-bound} hold under the assumption that $\psi(x)$ is $p$ times
continuously differentiable with Lipschitz continuous $p$-th derivative on the
``path of iterates'' $\cup_{k\geq 0}[x_k,x_{k+1}]$, in that
\beqn{segment-cond}
\bigmax_{\xi \in [0,1]} \| \nabla_x^p \psi(x_k+\xi s_k) - \nabla_x^p \psi(x_k) \|_p
\leq L_{\psi,p} \|s_k\|,
\eeqn
for all $\xi \in [0,1]$, all $k\in \calS$ and for some constant $L_{f,p} \geq
0$ independent of $x_k$ and $s_k$. (Obviously, if the $p$-th
derivative of $\psi$ is Lipschitz continuous in an open set containing $\calF$
or containing the level set $\{x \in \calF \mid \psi(x) \leq \psi(x_0)\}$,  
then \req{segment-cond} holds.)

At variance with the method proposed in \cite{CartGoulToin15b}, the algorithm
described in \cite{CartGoulToin16a} is of trust-region type with non-increasing
radius.  It approximately
minimizes a $q$-th degree Taylor inside such a region, Lemma~4.3 in
\cite{CartGoulToin16a} then ensures that, as long as \req{cc-approx} fails
(for general $q\geq 1$ this time and for $\Delta$ being the trust-region
radius at iteration $k$),
\beqn{cc-decr-succ-q}
\psi(x_k) - \psi(x_{k+1}) \geq \kap{decr}^\psi \epsilon^{q+1}
\eeqn
for each $k \in \calS$, where we have redefined the constant
$\kap{decr}^\psi$ to reflect the change in algorithm. In addition, Lemma~4.1
in the same paper also ensures that \req{cc-unsucc-bound} holds
for a redefined $\kap{uns}^\psi$. Both of these properties again hold if
$\psi(x)$ is $q$ times continuously differentiable with Lipschitz continuous
$q$-th derivative on the ``path of iterates'' $\cup_{k\geq 0}[x_k,x_{k+1}]$,
in the sense that \req{segment-cond} (with $p$ replaced by $q$). 

Summarizing, we see that there exist algorithms for the solution of
\req{ccprob} which use truncated Taylor series model of degree $q$ and ensure,
under suitable assumptions, both \req{cc-unsucc-bound} and, as long as
\req{cc-approx} does not hold for some algorithm-dependent non-increasing
$\Delta \in (0,1]$, a lower bound on the objective-function decrease at
successful iterations of the form
\beqn{cc-decr}
\psi(x_k) - \psi(x_{k+1}
\geq \kap{decr}^\psi[ \ed\s{LS} ]^\pi
\tim{ for }
k \in \calS
\eeqn
for suitable method-dependent constant $\kap{decr}^\psi \in (0,1)$ and parameter
$\pi\geq 1$. (We have that $\pi= (p+1)/p$ in \req{cc-decr-succ-1} and $\pi = q+1$ in
\req{cc-decr-succ-q}.)

Let us now turn to least-squares problems of the form
\beqn{lsprob}
\min_{x \in \calF} \psi(x) \eqdef \half \|F(x)\|^2,
\eeqn
(that is problem \req{ccprob} where $\psi(x) = \half \|F(x)\|_2^2$),
where $F$ is a smooth function from $\Re^n$ to $\Re^m$. Following \cite{CartGoulToin12e} and
\cite{CartGoulToin16a}, an $\epsilon$-approximate\footnote{$\ep\s{LS}$ is the 
primal accuracy for solving problem \req{lsprob} and $\ed\s{LS}$ the dual
one.} $q$-th order critical point for this problem can be defined as a point $x$
such that
\beqn{ls-approx}
\|F(x)\| \leq \ep\s{LS}
\tim{ or }
\phi_{\psi,j}^\Delta(x) \leq \ed\s{LS} \Delta^j \|F(x)\|
\tim{for}
j = 1,2
\eeqn
and some $\Delta \in (0,1]$.  Note that the second part of \req{ls-approx} has
the same form as \req{cc-approx} with $\epsilon$ in the former being replaced
by $\ed\s{LS} \|F(x)\|$ in the latter. As in
\cite{CartGoulToin12e,CartGoulToin15b}, it is now easy to verify that,
whenever $\|F(x_k)\| \geq \|F(x_{k+1})\|$ and as long as \req{ls-approx} fails
for $x_{k+1}$, 
\beqn{ls-decr}
\begin{array}{lcl}
\|F(x_k)\|\, \left( \|F(x_k)\| - \|F(x_{k+1})\| \right)
& \geq & \half \left( \|F(x_k)\| + \|F(x_{k+1})\| \right)
               \left( \|F(x_k)\| - \|F(x_{k+1})\| \right)\\*[1ex]
& \geq & \half\|F(x_k)\|^2 - \half \|F(x_{k+1})\|^2 \\*[1ex]
&   =  & \psi(x_k) - \psi(x_{k+1}) \\*[1ex]
& \geq & \kap{decr}^\psi[\, \ed\s{LS} \|F(x_{k+1})\|\,]^\pi,
\end{array}
\eeqn
where we used \req{cc-decr} with the form of the second part of
\req{ls-approx} to derive the last inequality. We will use this last 
formulation of the guaranteed decrease for least-squares problems as a key
piece of our evaluation complexity analysis, together with
\req{cc-unsucc-bound} which is needed because the algorithms under
consideration require one objective-function evaluation per iteration and one
evaluation of its derivatives per successful iteration.

\subsection{The outer algorithm}
  
The idea of the two-phase framework which we now introduce is to first apply
one of the least-squares algorithms discussed above or any other method with
similar guarantees), which we call \inneralg,to the problem
\beqn{feas-prob}
\min_{x\in \calF} \nu(x) \eqdef \half \|c(x)\|^2.
\eeqn 
(of the form \req{lsprob} with $\psi = \nu$) for finding (under suitably
adapted assumptions) an approximate feasible point, if possible. If one is found, 
\inneralg\ is then applied to approximately solve the problem 
\beqn{phase2-problem}
\min_{x \in \calF} \mu(x,t_k) 
\eqdef \half \|r(x,t_k)\|^2  
\eqdef \half \left\|\left( \begin{array}{c}
 c(x) \\ f(x) - t_k \end{array} \right) \right\|^2
\eeqn
(again of the form \req{lsprob} with $\psi = \mu$) for some monotonically
decreasing sequence of ``targets'' $t_k$ ($k=1,\ldots$). The resulting
algorithm is described \vpageref{outeralg}. Observe that the recomputations of
$\phi_{\mu,j}(x_{k+1},t_{k+1})$ $(j \in \ii{q})$ in Step~2.(b) do not require
re-eval\-uating $f(x_{k+1})$ or $c(x_{k+1})$ or any of their derivatives. 
\enlargethispage{8ex}
\balgo{outeralg}{{\sc{outer}}: a two-phase algorithm for constrained optimization}{
A starting point $x_{-1}$ 
and a criticality order $q \in \{1,2,3\}$ (for both the
feasibility phase and  the optimization phase) are given, as well as a
constant $\delta \in (0,1)$. The primal and dual tolerances   
$
0 < \ep < 1 \tim{ and } 0 < \ed < 1
$
are also given.}{
\begin{description}
\item[Phase 1: ] \mbox{}\\
Starting from $x_0=P_\calF(x_{-1})$, apply \inneralg\ to minimize
$\nu(x)=\half \|c(x)\|^2$ subject to $x \in \calF$ until a point $x_1 \in \calF$
and $\Delta_0 \in (0,1]$ are found such that
\beqn{phase1-term}
\|c(x_1)\| < \delta \ep \tim{ or }
\phi_{\nu,j}^{\Delta_1}(x_1) \leq \ed \Delta_0^j \|c(x_1)\|
\ms
(j \in \ii{q}).
\eeqn
If $\|c(x_1)\| > \delta \ep$, terminate with $x_\epsilon = x_1$.

\item [Phase 2: ] \mbox{}
\begin{enumerate}
\item Set $t_1 = f(x_1)-\sqrt{\ep^2-\|c(x_1)\|^2}$.

\item For $k=1, 2, \ldots$, do:
\begin{enumerate}
\item Starting from $x_k$, apply \inneralg\ to minimize $\mu(x,t_k)$ as a
  function of $x \in \calF$ until an iterate $x_{k+1} \in \calF$ and $\Delta_k
  \in (0,\Delta_{k-1}]$ are found such that 
    \beqn{ph2term}
   \begin{array}{c}
   \|r(x_{k+1},t_k)\| <  \delta \ep
   \tim{ or }
    f(x_{k+1}) < t_k \\*[2ex]
    \tim{ or  }
    \phi_{\mu,j}^{\Delta_k}(x_{k+1},t_k) \leq \ed  \Delta_k^j \|r(x_{k+1},t_k)\|
    \ms (j \in \ii{q}).
    \end{array}
  \eeqn
\item
\begin{enumerate}
\item If $\|r(x_{k+1},t_k)\| <  \delta \ep$,
define $t_{k+1}$ according to 
\beqn{tk-update}
t_{k+1}=f(x_{k+1})-\sqrt{\ep^2-\|c(x_{k+1})\|^2}.
\eeqn
and terminate with $(x_\epsilon,t_\epsilon) = (x_{k+1},t_{k+1})$ if 
\beqn{ph2-termp}
\phi_{\mu,j}^{\Delta_k}(x_{k+1},t_{k+1}) \leq \ed \Delta_k^j \|r(x_{k+1},t_{k+1})\|
\tim{for} j \in \ii{q}.
\eeqn
\item If  $\|r(x_{k+1},t_k)\| \geq \delta \ep$ and  $f(x_{k+1}) < t_k$, 
define $t_{k+1}$ according to
\beqn{tk-swap}
t_{k+1} = 2 f(x_{k+1})-t_k
\eeqn
and terminate with $(x_\epsilon,t_\epsilon) = (x_{k+1},t_{k+1})$ if
\req{ph2-termp} holds.
\item If $\|r(x_{k+1},t_k)\| \geq \delta \ep$ and $f(x_{k+1}) \geq  t_k$,
terminate with $(x_\epsilon,t_\epsilon) = (x_{k+1},t_k)$
\end{enumerate}
\end{enumerate}
\end{enumerate}
\end{description}
}

We now derive some useful properties of \outeralg .  For this purpose,
we partition the Phase~2 outer iterations (before that where termination
occurs) into two subsets whose indexes are given by   
\beqn{K-pos}
\calK_+ 
\eqdef \{ k \geq 0 \mid \|r(x_{k+1},t_k)\| <  \delta \ep
          \tim{ and \req{tk-update} is applied} \}
\eeqn
and
\beqn{K-neg}
\calK_- 
\eqdef \{ k \geq 0 \mid \|r(x_{k+1},t_k)\| \geq  \delta \ep
          \tim{ and \req{tk-swap} is applied} \}
\eeqn
The partition \req{K-pos}-\req{K-neg} allows us to prove then following technical results.

\llem{tech-l}{
The sequence  $\{t_k\}$ is monotonically decreasing.  Moreover, in every Phase
2 iteration of \outeralg\ of index $k \geq 1$, we have that 
\beqn{tkltfk}
f(x_k) - t_k \geq 0,
\eeqn
\beqn{nreseps}
\|r(x_{k+1}, t_{k+1})\| = \ep  \tim{ for } k \in \calK_+,
\eeqn
\beqn{r-swap}
\| r(x_{k+1},t_{k+1})\| = \| r(x_{k+1},t_k)\| \leq \ep \tim{ for } k \in \calK_-,
\eeqn
\beqn{approxfeas}
\|c(x_k)\| \leq \ep \tim{and} f(x_k) - t_k \leq \ep,
\eeqn 
\beqn{tkdecr}
t_k - t_{k+1} \geq (1-\delta) \ep \tim{ for } k \in \calK_+.
\eeqn
Finally, at termination of \outeralg,
\beqn{conds-at-termination}
\begin{array}{c}
\|r(x_\epsilon,t_\epsilon)\| \geq  \delta \ep,
\ms f(x_\epsilon) \geq t_\epsilon \\*[2ex]
\tim{ and }
\phi_{\mu,j}^{\Delta_k}(x_\epsilon,t_\epsilon) \leq\ed \Delta_k^q \|r(x_\epsilon,t_\epsilon)\|
\tim{for} j \in \ii{q}.
\end{array}
\eeqn
}

\proof{
The inequality \req{tkltfk} follows from \req{tk-update} for $k-1 \in \calK_+$
and from \req{tk-swap} for $k-1 \in \calK_-$.  \req{nreseps} is also deduced
from \req{tk-update} while  \req{tk-swap} implies the equality in
\req{r-swap}, the inequality in that statement resulting from the
monotonically decreasing nature of $\|r(x,t_k)\|$ during inner iterations in
Step~2.(a) of \outeralg. The inequalities \req{approxfeas} then
follow from \req{tkltfk},  \req{nreseps} and \req{r-swap}. 
We now prove \req{tkdecr}, which only occurs when  $\|r(x_{k+1},t_k)\| \leq
\delta \ep$, that is when 
\beqn{res1}
(f(x_{k+1}) - t_k)^2 + \|c(x_{k+1})\|^2 \leq \delta^2 \ep^2.
\eeqn
From \req{tk-update}, we then have that
\beqn{dtk}
t_k-t_{k+1}= -(f(x_{k+1})-t_k)+\sqrt{\|r(x_k,t_k)\|^2-\|c(x_{k+1})\|^2}.
\eeqn
Now taking into account that the global minimum of the problem 
\[
\min_{(f,c)\in\smallRe^2} \vartheta(f,c)\eqdef -f+\sqrt{\ep^2 - c^2} 
\tim{subject to} f^2+c^2 \leq \omega^2,
\]
for $\omega \in [0, \ep]$ is attained at $(f_*,c_*)=(\omega, 0)$ and it is given by
$\vartheta(f_*,c_*)=\ep-\omega$ (see \cite[Lemma 5.2]{CartGoulToin13a}), we obtain
from \req{res1} and \req{dtk} (setting $\omega = \delta \ep$) that
\[
t_k-t_{k+1} \geq \ep - \omega = (1-\delta)\ep \tim{ for } k \in \calK_+
\]
for $k \in \calK_+$, which is \req{tkdecr}.
Note that, if $k \in \calK_-$, then we must have that $t_k > f(x_{k+1})$ and
thus \req{tk-swap} ensures that $t_{k+1} < t_k$. This observation and
\req{tkdecr} then allow us to conclude that the sequence $\{t_k\}$ is
monotonically decreasing. 

In order to prove \req{conds-at-termination}, we need to consider, in turn,
each of the three possible cases where termination occurs in Step~2.(b).  In
the first case (i), $\|r(x_{k+1},t_k)\|$ is small (in the sense that the first
inequality in \req{ph2term} holds) and \req{tk-update} is then used, implying
that \req{nreseps} holds and that $f(x_{k+1}) > t_{k+1}$. If termination
occurs because \req{ph2-termp} holds,
then \req{conds-at-termination} clearly holds at $(x_{k+1},t_{k+1})$.  In the
second case (ii), the residual $\|r(x_{k+1},t_k)\|$ is large (the first
inequality in \req{ph2term} fails), but $f(x_{k+1}) < t_k$, and $t_{k+1}$ is
then defined by \req{tk-swap}, ensuring that $f(x_{k+1}) > t_{k+1}$ and,
because of \req{r-swap}, that $\|r(x_{k+1},t_{k+1})\|$ is also large.  As
before \req{conds-at-termination} holds at $(x_{k+1},t_{k+1})$ if termination
occurs because \req{ph2-termp} is satisfied.
The third case (iii) is when
$\|r(x_{k+1},t_k)\|$ is sufficiently large and $f(x_{k+1}) \geq t_k$.  But
\req{ph2term} then guarantees that $\phi_{\mu,j}^{\Delta_k}(x_{k+1},t_k) \leq \ed
\Delta_k^j\|r(x_{k+1},t_k)\|$ for $j \in \ii{q}$, and the inequalities
\req{conds-at-termination} are again satisfied at $(x_{k+1},t_k)$.
} 

\numsection{Evaluation complexity}\label{complexity-s}

In order to state the smoothness assumptions for problem \req{genprob}, we
first define, for some parameter $\beta >0$, the  neighbourhood of the
feasible set given  by 
\[
\calC_\beta= \{ x \in \calF \mid \|c(x)\| \leq \beta \}.
\]
We then assume the following.

\ass{AS.1}{ The feasible set $\calF$ is closed, convex and non-empty.}

\ass{AS.2}{The function $\nu(x)$ is smooth enough to ensure that conditions
  \req{ls-decr} and \req{cc-unsucc-bound} hold for \inneralg\ applied on
  problem \req{feas-prob}.
}
\ass{AS.3}{The function $\mu(x,t)$ is smooth enough in $x$ to ensure that conditions
  \req{ls-decr} and \req{cc-unsucc-bound} hold for \inneralg\ applied on
  problem \req{phase2-problem}, with constants $\kap{decr}^\mu$ and
  $\kap{uns}^\mu$ independent of $t$.
}

\ass{AS.4}{There exists constants $\beta \geq \ep$ and $f_{\rm low}\in \Re$
  such that $f(x)\geq f_{\rm low}$ for all $x \in \calC_\beta \eqdef \{ x \in
  \calF \mid \, \|c(x)\| \leq \beta\}$.
}

\noindent
AS.2 and AS.3 remain implicit and depend on the particular inner algorithm
used (see Section~\ref{inner-s}).  For completeness, we now give conditions on
the problem's functions $f$ and $\{c_i\}_{i=1}^m$ which allow the transition between
assumptions on $f$ and $c$ and the required ones on the Phase~1 and Phase~2
objective functions $\nu$ and $\mu$.

\llem{LS-Lip}{Let $p\geq 1$. Assume that $f$ and $\{c_i\}_{i=1}^m$ are $p$
  times continuoulsy differentiable and that their derivatives of order one up
  to $p$ are uniformly bounded and Lipschitz continuous in an open set
  containing $\calF$. Let the iterations of \inneralg\ applied to problem
  \req{feas-prob} be indexed by $j$. Then \req{segment-cond} holds for
  $\nabla_x^q\nu(x)$ on every segment $[x_j, x_j+s_j]$ $(j\geq 0)$ generated
  by \inneralg\ during Phase~1 and any $q\in\ii{p}$. The same conclusion holds
  for $\nabla_x^q \mu(x,t)$ on every segment $[x_j, x_j+s_j]$ $(j\geq 0)$
  generated by \inneralg\ during Step~2.(a) of Phase~2 and any $q\in\ii{p}$,
  the Lipschitz constant in this latter case being independent of $t$.
}

\proof{
Since
\[
\nabla_x^q \nu(x) = \sum_{i=1}^m \left[
  \sum_{\ell,j>0, \, \ell+j=q} \alpha_{\ell,j} \nabla_x^j c_i(x) \otimes \nabla_x^\ell c_i(x)
  + c_i(x) \nabla_x^q c_i(x) \right]
\]
(where $\{\alpha_{\ell,j}\}$ are suitable non-negative and finite coefficients),
condition \req{segment-cond} is satisfied on the segment $[x_j, x_j+s_j]$ if
(i) the derivatives $\{\nabla_x^{\min[\ell,j]}c_i(x)\}_{i=1}^m$ are Lipschitz continuous
    on $[x_j, x_j+s_j]$,
(ii) $\{\nabla_x^{\max[\ell,j]}c_i(x)\}_{i=1}^m$ are uniformly bounded
    on $[x_j, x_j+s_j]$, and
(iii) we have that
    \beqn{cnablaq}
    \sum_{i=1}^m \|c_i(x_j+\xi s_j)
                \nabla_x^qc_i(x_j+\xi s_j)-c_i(x_j)\nabla_x^qc_i(x_j)\|_q
    \leq L_1 \xi \|s_j\|
    \eeqn
    for some constant $L_1 > 0$.
    The first two of these conditions are ensured by the lemma's assumptions.
Moreover,
\[
\begin{array}{lcl}
  \lefteqn{\|c_i(x_j+\xi s_j)\nabla_x^q c_i(x_j+\xi s_j)
             - c_i(x_j)\nabla_x^q c_i(x_j) \|_q} \\
\ms\ms\ms\ms\ms\ms\ms\ms
& \leq & |c_i(x_j+\xi s_j)-c_i(x_j)|\,\|\nabla_x^q c_i(x_j+\xi s_j)\|_q \\
&      &  + |c_i(x_j)| \,\|\nabla_x^q c_i(x_j+\xi s_j)-\nabla_x^q c_i(x_j)\|_q
\end{array}
\]
and the first term on the right-hand side is bounded above by $L^2\xi\|s_j\|$
and the second by $|c_i(x_j)|L\xi\|s_j\|$.  Hence \req{cnablaq} holds with
\[
L_1 = \sum_{i=1}^m \left(L^2+|c_i(x_j)|L\right)
\leq mL^2 + m\|c(x_j)\|L
\leq mL^2 + m\|c(x_0)\|L
\]
because \inneralg\  ensures that $\|c(x_j)\|\leq \|c(x_0)\|$ for
all $j\geq 0$. As a consequence, the lemma's assumptions guarantee that
\req{segment-cond} holds with the Lipschitz constant
\[
m\left[\left(\max_{i=1,\ldots,m}\alpha_i\right) L^2 + L^2 + \|c(x_0)\|L\right].
\]
We may now repeat, for $\mu(x,t)$
(with fixed $t$) the same reasoning as above and obtain that condition
\req{segment-cond} holds for each segment
$[x_j,x_j+s_j]$ generated by \inneralg\ applied in Step≃2.(a) of
Phase~2, with Lipschitz constant 
\[
\begin{array}{lcl}
\lefteqn{m\left[\left(\bigmax_{i=1,\ldots,m}\alpha_i \right)L^2 + L^2
  + \|c(x_{j,0})\| L \right] +\left(\bigmax_{i=1,\ldots,m}\alpha_i \right)L^2
  + L^2 + |f(x_{j,0})-t_j|L }\\*[2ex]
\ms\ms\ms\ms\ms\ms\ms\ms\ms\ms
& \leq & \!\! (m+1)\left[L^2 \left(1+ \bigmax_{i=1,\ldots,m}\alpha_i
  \right) + L \right]
\eqdef L_{\mu,p},
\end{array}
\]
where we have used \req{nreseps} and $\ep\leq 1$ to deduce the inequality.
Note that this constant is independent of $t_j$, as requested.
}

\noindent
As the constants $\kap{decr}^\mu$ and $\kap{uns}^\mu$ in \req{ls-decr} and
\req{cc-unsucc-bound} directly depend, for the class of inner algorithms
considered, on the Lipschitz constants of the derivatives of $\mu$ with
respect to $x$, the independence of these with respect to $t$ ensures that
$\kap{decr}^\mu$ and $\kap{uns}^\mu$ are also independent of $t$, as requested
in AS.3.

We now start the evaluation complexity analysis by examining the complexity of
Phase~1 of \outeralg.

\llem{ph1-compl}{
Suppose that AS.1 and AS.2 hold.  Then Phase 1 of
\outeralg\ terminates with an $x_1$ such that $\|c(x_1)\|\leq \delta
\ep$ or $\phi_{\nu,q}^{\Delta_k} \leq \epsilon \Delta_k^q$ after at most 
\[
\left\lfloor
\kap{CC}^{\|c\|} \|c(x_0)\| \, \max\left[\ep^{-1},
                                        \ep^{1-\pi} \, \ed^{-\pi}
                                 \right] \right\rfloor +1
\]
evaluations of $c$ and its derivatives, where
$
\kap{CC}^{\|c\|} \eqdef 2^{-\pi} \kappa_u[\kap{decr}^\nu]^{-1} \delta^{1-\pi}
$
with $\kap{decr}^\nu$ being the problem-dependent constant defined in
\req{ls-decr} for the function $\nu(x)$ corresponding to \req{feas-prob}. 
}

\proof{First observe that, as long as \inneralg\ applied to
problem \req{feas-prob} has not terminated,
\beqn{biggish}
\|c(x_\ell)\| \geq \delta \ep,
\eeqn
because of the first part of \req{phase1-term}.
Let $\ell \in \calS_k$ be the index of a successful iteration of
\inneralg\ before termination and suppose first that $\|c(x_{\ell+1})\| \leq
\half \|c(x_\ell)\|$.  Then
\beqn{bcj1}
\|c(x_\ell)\| - \|c(x_{\ell+1})\| \geq \half \|c(x_\ell)\| \geq \half \delta \,\ep
\eeqn
Suppose now that $\|c(x_{\ell+1})\| > \half \|c(x_\ell)\|$.
As a consequence, we obtain that
\[
(\|c(x_\ell)\| - \|c(x_{\ell+1})\|)\,\|c(x_\ell)\|
\geq  \kap{decr}^{\nu}\, (\ed \|c(x_{\ell+1})\|)^{\pi}
\]
where we have also the fact that
$\phi_{\nu,j}^{\Delta_k}(x_{\ell+1}) > \ed\|c(x_{\ell+1})\|\Delta_k^j$
since $\ell$ occurs before termination, the fact that $\|c(x_\ell)\| \geq
\|c(x_{\ell+1})\|$ for $\ell \in \calS$ and condition \req{ls-decr}.  Hence,
using \req{biggish}, we have that
\[
\|c(x_\ell)\| - \|c(x_{\ell+1})\|
\geq \kap{decr}^\nu  2^{-\pi} \|c(x_\ell)\|^{\pi-1} \, \ed^{\pi}
\geq 2^{-\pi} \kap{decr}^\nu\,\delta^{\pi-1}\;\ep^{\pi-1}\,\ed^{\pi}.
\]
Because of the definition of $\kap{decr}^\nu$ in \req{ls-decr}, we thus obtain
from this last bound and \req{bcj1} that,  for all $j$,  
\[
\|c(x_\ell)\| - \|c(x_{\ell+1})\|
\geq \half\kap{decr}^\nu\,\delta^{\pi-1}
\min\left[ \ep, \ep^{\pi-1}\,\ed^{\pi} \right].
\]
We then deduce that
\[
|\calS_k|
\leq 2[\kap{decr}^\nu]^{-1}\delta^{-\frac{1}{p}} \, \|c(x_0)\| \,
\max\left[\ep^{-1},\ep^{1-\pi}\,\ed^{-\pi}\right]
\]
The desired conclusion then follows by using condition \req{cc-unsucc-bound}
and adding one for the final evaluation at termination.
} 

\noindent
Using the results of this lemma allows us to bound the number of outer
iterations in $\calK_+$.

\llem{number-outer}{
Suppose that AS.4 holds. Then 
\[
| \calK_+ | \leq \frac{f(x_1)-f_{\rm low}+1}{1-\delta}\,\ep^{-1}.
\]
}

\proof{
We first note that \req{nreseps} and \req{r-swap} and AS.4 ensure that $x_k \in
\calC_\beta$ for all $k \geq 0$. The result then immediately follows from AS.4
again and the observation that, from \req{tkdecr}, $t_k$ decreases monotonically
with a decrease of at least $(1-\delta)\ep$ for $k \in \calK_+$.
} 

\noindent
Consider now $x_k$ for $k \in \calK_+$ and  denote by $x_{n(k)}$ the next
iterate such that $n(k) \in \calK_+$ or the algorithm terminates at
$n(k)$. Two cases are then possible: either a single pass
in Step~2.(a) of \outeralg\ is sufficient to obtain
$x_{n(k)}$ ($n(k) = k+1$) or  two or more passes are necessary, with
iterations $k+1, \ldots, n(k)-1$ belonging to $\calK_-$.  Assume now that
the iterations of \inneralg\ at Step~2.(a) of the outer
iteration $\ell$ are numbered $(\ell,0), \,(\ell,1), \ldots, (\ell, e_\ell)$
and note that the mechanism of \outeralg\ ensures that iteration
$(\ell,e_\ell)$ is successful for all $\ell$.  Now define, for $k \in
\calK_+$, the index set of all inner iterations necessary to deduce $x_{n(k)}$
from $x_k$, that is 
\beqn{Ikdef}
\calI_k
\eqdef \{ (k,0), \ldots, (k,e_k), \ldots, (\ell,0), \ldots, (\ell,e_\ell),
          \ldots, (n(k)-1,0), \ldots (n(k)-1,e_{n(k)-1}) \}
\eeqn
where $k < \ell < n(k)-1$.
Observe that, by the definitions  \req{K-pos} and \req{Ikdef}, the index set
of all inner iterations before termination is given by $\cup_{k \in \calK_+}
\calI_k$, and therefore that the number of evaluations of problem's functions
required to terminate in Phase~2 is bounded above by  
\beqn{rough-bound}
|\bigcup_{k \in \calK_+} \calI_k| +1
\leq \Bigg(\frac{f(x_1)-f_{\rm low}+1}
                {1-\delta}\ep^{-1} \times \max_{k \in \calK_+}|\calI_k| \Bigg)+1,
\eeqn
where we added 1 to take the final evaluation into account and where we used
Lemma~\ref{number-outer} to deduce the inequality.  We now invoke the
complexity properties of \inneralg\ applied to problem
\req{phase2-problem} to obtain an upper bound on the cardinality of each
$\calI_k$.

\llem{number-inner}{
  Suppose that AS.1--AS.3 hold. Then, for each $k \in \calK_+$ before
  termination, 
\[
|\calI_k|
\leq (1-\delta)\kap{CC}^\mu\,\max\left[1, \ep^{2-\pi}\ed^{-\pi}\right].
\]
where $\kap{CC}^\mu$ is independent of $\ep$ and $\ed$ and captures the
problem-dependent constants associated with problem \req{phase2-problem} for
all values of $t_k$ generated by the algorithm.
}

\proof{
Observe that \req{r-swap} and the
mechanism of this algorithm guarantees the strictly decreasing nature of the
sequence $\{\|r(x_\ell,t_\ell)\|\}_{\ell = k}^{n(k)-1}$ and hence of the
sequence $\{\|r(x_{\ell,s},t_\ell)\|\}_{(\ell,s) \in \calI_k}$. For each $k
\in \calK_+$, this reduction starts from the initial value $\|r(x_{k,0},t_k)\|=
\ep$ and is carried out for all iterations with index in $\calI_k$ at
worst until it is smaller than $\delta\ep$ (see the first part of
\req{ph2term}) or $\phi_{\mu,j}(x_{\ell,s}) \leq \ed
\Delta_k^j\|r(x_{\ell,s+1},t_\ell)\|$ for $j \in \ii{q}$.
We may then invoke \req{phase2-problem} and
\req{ls-decr} to deduce that, if
$(k,s)\in \calI_k$,
\beqn{decrrj}
(\|r(x_{k,s},t_k)\| - \|r(x_{k,s+1},t_k)\|) \|r(x_{k,s},t_k)\|
\geq \kap{decr}^\mu (\ed \|r(x_{k,s+1},t_k)\|)^{\pi},
\eeqn
for $0 \leq s<e_k$, while
\[
\half \|r(x_{k,e_k},t_k)\| - \half \|r(x_{k+1,0},t_{k+1})\| = 0.
\]
As above, suppose first that
$\|r(x_{k,s+1},t_k)\| \leq \half\|r(x_{k,s},t_k)\|$. Then
\beqn{brj1}
\|r(x_{k,s},t_k)\| - \|r(x_{k,s+1},t_k)\|
\geq \half \|r(x_{k,s},t_kl)\|
\geq \half \delta \ep
\eeqn
because of the first part of \req{ph2term}. If $\|r(x_{k,s+1},t_k)\| >
\half\|r(x_{k,s},t_k)\|$ instead, then \req{decrrj} implies that
\[
\|r(x_{k,s},t_k)\| - \|r(x_{k,s+1},t_k)\|
\geq  \kap{decr}^\mu \,2^{-\pi}\|r(x_{k,s},t_k)\|^{\pi-1} \, \ed^{\pi}.
\]
Combining this bound with \req{brj1} gives that 
\[
\|r(x_{k,s},t_k)\| - \|r(x_{k,s+1},t_k)\|
\geq 2^{-\pi} \kap{decr}^\mu \delta^{\pi-1}\;
\min\left[\ep, \ep^{\pi-1} \ed^{\pi}\right].
\]
and therefore, as in Lemma~\ref{ph1-compl}, that
\[
|\calI_k|
\leq 2^\pi [\kap{decr}^\mu]^{-1} \delta^{1-\pi}
\left[\frac{\ep -\delta \ep}
           {\min\left[\ep, \ep^{\pi-1} \ed^{\pi}\right]}\right]
= 2^\pi (1-\delta) \delta^{1-\pi}[\kap{decr}^\mu]^{-1}\,
       \max\left[1, \ep^{2-\pi}\ed^{-\pi}\right],
\]
and the conclusion follows with $\kap{CC}^\mu\eqdef 2^\pi \delta^{1-\pi}[\kap{decr}^\mu]^{-1}$.
} 

We finally combine the above results in a final theorem stating an evaluation
complexity bound for \outeralg\ in terms of the measures
$\phi_{\nu,j}^{\Delta_k}(x_\epsilon)$.

\lthm{general-th}{
Suppose that AS.1--AS.4 hold. Then, for some
constants $\kap{CC}^{\|c\|}$ and $\kap{CC}^\mu$ independent of $\ep$ and
$\ed$, \outeralg\   applied to problem \req{genprob} needs at most 
\beqn{finbound}
\left\lfloor
\Bigg(\kap{CC}^{\|c\|} \|c(x_0)\| + \kap{CC}^\mu [f(x_1)-f_{\rm low}+1] \Bigg)\,
\max\left[\ep^{-1}, \ep^{1-\pi}\ed^{-\pi}\right]
\right\rfloor + 2
\eeqn
evaluations of $f$, $c$ and their derivatives up to order $p$ to compute a
point $x_\epsilon$ and (possibly) a $t_\epsilon \leq f(x_\epsilon)$ such that,
when $t_\epsilon = f(x_\epsilon)$, 
\beqn{xinfeas}
\|c(x_\epsilon)\| >\delta \ep,
\tim{ and }
\phi_{\nu,j}^{\Delta_k}(x_\epsilon) \leq \ed\Delta_k^j\|c(x_\epsilon)\|
\tim{for} j \in \ii{q}
\eeqn
or, when $t_\epsilon < f(x_\epsilon)$,
\beqn{mustat}
\|c(x_\epsilon)\| \leq \ep,
\tim{ and }
\phi_{\mu,j}^{\Delta_k}(x_\epsilon,t_\epsilon) \leq \ed \Delta_k^j \|r(x_\epsilon,t_\epsilon)\|
\tim{for} j \in \ii{q}.
\eeqn
}

\proof{
If \outeralg\ terminates in Phase~1, we immediately obtain that
\req{xinfeas} holds, and Lemma~\ref{ph1-compl} then ensures that the number of
evaluations of $c$ and its derivatives cannot exceed 
\beqn{bph1}
\left\lfloor \kap{CC}^{\|c\|} \|c(x_0)\| \,
   \max\left[\ep^{-1}, \ep^{1-\pi}\ed^{-\pi}\right] \right\rfloor
+1.
\eeqn
The conclusions of the theorem therefore hold in this case.
Let us now assume that termination does not occur in Phase~1.  Then 
\outeralg\ must terminate after a number of evaluations of $f$ and
$c$ and their derivatives which is bounded above by the upper bound on the
number of evaluations in Phase~1 given by \req{bph1}  plus the bound on the
number of evaluations of $\mu$  given by \req{rough-bound} and
Lemma~\ref{number-inner}. Using the inequality $q_\nu \leq q$ and the facts
that  $\lfloor a \rfloor + \lfloor b \rfloor \leq \lfloor a+b \rfloor$ for
$a,b \geq 0$ and $\lfloor a +i \rfloor = \lfloor a \rfloor + i$ for $a \geq 0$
and $i \in \Na$, this yields the combined upper bound  
\[
\begin{array}{l}
\left\lfloor
  \kap{CC}^{\|c\|} \|c(x_0)\|\,
  \max \left[ \ep^{-1}, \ep^{1-\pi}\ed^{-\pi} \right]
 \right.
  \\*[2ex]
  \ms \ms \ms \ms  + \left.\left[(1-\delta)\kap{CC}^\mu 
    \,\max\left[1, \ep^{2-\pi} \ed^{-\pi}\right]\right]\times
  \left[ \bigfrac{f(x_1)-f_{\rm low}+1}{1-\delta} \,
          \ep^{-1}\right] \right\rfloor + 2,
\end{array}
\]
and \req{finbound} follows.
Remember now that \req{conds-at-termination} holds at termination of Phase~2,
and therefore that  
\beqn{rbiggish}
\ep \geq \|r(x_\epsilon,t_\epsilon)\| \geq  \delta \ep.
\eeqn
Moreover,  we also obtain from \req{conds-at-termination} that
\beqn{chi-small}
\phi_{\mu,j}^{\Delta_k}(x_\epsilon,t_\epsilon) 
\leq \ed \Delta_k^j \|r(x_\epsilon,t_\epsilon)\|
\tim{ for } j \in \ii{q}.
\eeqn
Assume first that $f(x_\epsilon) = t_\epsilon$.  Then, using  the definition
of $r(x,t)$, we deduce that, for $j \in \ii{q}$,
\[
\phi_{\nu,j}^{\Delta_k}(x_\epsilon) 
= \phi_{\mu,j}^{\Delta_k}(x_\epsilon) 
\leq \ed \Delta_k^j \|c(x_\epsilon)\|
\]
and \req{xinfeas} is again satisfied because \req{rbiggish} gives
that
$\|c(x_\epsilon)\| = \|r(x_\epsilon,t_\epsilon)\| \geq  \delta \ep$.

If $f(x_\epsilon) > t_\epsilon$ (the case where $f(x_\epsilon) <
t_\epsilon$ is excluded by \req{conds-at-termination}), we see that 
the inequality $\|c(x_\epsilon)\| \leq \|r(x_\epsilon,t_\epsilon)\| \leq \ep$,
and \req{chi-small} imply \req{mustat}.
} 

\noindent
Note that the bound \req{finbound} is
$O(\epsilon^{-(2\pi-1)})$ whenever $\ep = \ed = \epsilon$.
Also note that we have used the same algorithm for Phase~1 and Phase~2 of
\outeralg, but we could choose to use different methods of complexity
$\pi_\nu$ and $\pi_\mu$, respectively, leading a final bound of the form
\beqn{twophase-bound}
O\left(\max\left[\ep^{-1}, \ep^{1-\pi_\nu}\ed^{-\pi_\nu}\right]
+ \max\left[\ep^{-1}, \ep^{1-\pi_\mu}\ed^{-\pi_\mu}\right] \right).
\eeqn
Different criticality order may also be chosen for the two phases, leading to
variety of possible complexity outcomes.

It is important to note that the complexity bound given by
Theorem~\ref{general-th} depends linearly on $f(x_1)$, the value of the
objective function at the end of Phase~1.  Giving an $\epsilon$-independent
upper bound on this quantity is in general impossible, but can be done in some
case. A trivial bound can of course be obtained if $f(x)$ is bounded in a
neighbourhood of the feasible set, that is $\{ x \in \calF | \|c(x)\| \leq
\beta \}$ for some $\beta >0$. This has the advantage of providing a
complexity result which is self-contained (in that it only involves
problem-dependent quantities), but it is quite restrictive as it excludes, for
instance, problems with equality constraints only ($\calF = \Re^n$) and
coercive objective functions. A bound is also readily obtained if the set
$\calF$ is itself bounded (for instance when the variables are subject to
finite lower and upper bounds) or if one assumes that the iterates generated
by Phase~1 remain bounded.  This may for example be the case if the set $\{ x
\in \Re^n \mid c(x) = 0 \}$ is bounded. For specific choices of the
convexly-constrained algorithm applied for Phase~1 of \outeralg, an
$\ep$-dependent bound can finally be obtained without any further
assumption. If Phase~1 is solved using the trust-region based algorithm of
\cite{CartGoulToin16a} and $x_1$ is produced after $k_\epsilon$ iterations of
this algorithm, we obtain from the definition of the step that $\|s_k\| \leq
\Delta_{\max}$ for all $k \geq 1$. In the same spirit, if the regularization
algorithm of \cite{CartGoulToin15b} is used for Phase~1 and $x_1$ is produced
after $k_\epsilon$ iterations of this algorithm, we obtain from the proof of
Lemma~2.4 in \cite{CartGoulToin15b} and the definition of successful iterations that  
\[
\nu(x_0) \geq \nu(x_0)-\nu(x_1)
= \sum_{k\in \calS_{k_\epsilon}}[\nu(x_k)-\nu(x_{k+1})]
\geq \frac{\eta \sigma_{\min}}{(p+1)!}\sum_{k\in \calS_{k_\epsilon}}\|s_k\|^{p+1},
\]
giving that
\[
\|s_k\|\leq \left(\frac{\nu(x_0)(p+1)!}{\eta  \sigma_{\min}}\right)^{\frac{1}{p+1}}.
\]
Hence $\|x_1-x_0\|$ is itself bounded above by this
constant times the ($\ep$-dependent) number of iterations in Phase~1 given
by Lemma~\ref{ph1-compl}. Using the boundedness of the gradient of $\nu(x)$
on the path of successful iterates implied by AS.2 then ensures (see
Appendix) the (extremely pessimistic) upper bound 
\beqn{fx1-bound}
f(x_1)=f(x_0)+O\left(\max\left[\ep^{-1},\ep^{1-\pi}\,\ed^{-\pi}\right]\right).
\eeqn
Substituting this bound in
\req{finbound} in effect squares the complexity of obtaining
$(x_\epsilon,t_\epsilon)$.

Assuming that $f(x_1)-f_{\rm low}$ can be bounded by a constant independent of $\ep$
and $\ed$, Table~\ref{compl-table} gives the evaluation complexity bound for achieving
first-and second-order optimality for the problem with additional equality constraints,
depending on the choice of underlying algorithm for convexly-constrained
optimization. In this table, $q$ is the sought criticality order and $p$ is the
degree of the Taylor series being used to model the objective function in the
inner algorithm. The table also shows that the use of regularized high-degree
models for optimality orders beyond one remains to be explored.

\bctable{c|l|ccc}
      & \multicolumn{1}{|c}{TR-algo}  & \multicolumn{3}{|c}{Regularization} \\
 $q$  & \multicolumn{1}{|c}{$(p=q)$}  & \multicolumn{1}{|c}{$p=q$} &
        \multicolumn{1}{c}{$p=q+1$}   & \multicolumn{1}{c}{$p \geq  q$} \\
\hline
  1   & $O\Big(\epsilon^{-3}\Big)$     &  $O\Big(\epsilon^{-3}\Big)$  &
        $O\Big(\epsilon^{-2}\Big)$     &  $O\Big(\epsilon^{-\frac{p+2}{p}}\Big)$ \\
  2   & $O\Big(\epsilon^{-5}\Big)$     &  ? & ? & ? \\
 $q$  & $O\Big(\epsilon^{-(2q+1)}\Big)$ &  ? & ? & ?
\ectable{\label{compl-table}Evaluation complexity bounds for \outeralg\ as a
  function of the underlying algorithm for convexly-constrained problems, for
  $\epsilon$-independent $f(x_1)-f_{\rm low}$ and $\epsilon = \ep = \ed$}

We now consider the link between the necessary conditions derived in
Section~\ref{cnc-s} and the results of Theorem~\ref{general-th}. For future reference,
we start by giving the full expressions of the first four derivatives of
$\mu(x,t)$ as a function of $x$:
\beqn{D1M}
\nabla_x^1 \mu(x,t) = \sum_{i=1}^m c_i(x) \nabla_x^1c_i(x) + (f(x)-t)\nabla_x^1f(x),
\eeqn
\beqn{D2M}
\nabla_x^2 \mu(x,t) = \sum_{i=1}^m \Big[\nabla_x^1c_i(x)\otimes \nabla_x^1c_i(x)
  + c_i(x)\nabla_x^2c_i(x)\Big] + \nabla_x^1f(x)\otimes \nabla_x^1f(x)
  + (f(x)-t)\nabla_x^2f(x)
\eeqn
\beqn{D3M}
\nabla_x^3 \mu(x,t)
= \sum_{i=1}^m \Big[ 3\,\nabla_x^2c_i(x)\otimes\nabla_x^1c_i(x)
   +c_i(x)\nabla_x^3c_i(x)\Big] 
   +  3\,\nabla_x^2f(x)\otimes\nabla_x^1f(x)+(f(x)-t)\nabla_x^3f_(x)
   \eeqn
\beqn{D4M}
\begin{array}{ll}
\lefteqn{ \nabla_x^4 \mu(x,t)
= \sum_{i=1}^m \Big[4\,\nabla_x^3c_i(x)\otimes\nabla_x^1c_i(x)
  +3\,\nabla_x^2c_i(x)\otimes\nabla_x^2c_i(x) +c_i(x)\nabla_x^4c_i(x)\Big]}\\
\ms\ms\ms\ms\ms\ms\ms\ms\ms\ms
& +4\,\nabla_x^3f(x)\otimes\nabla_x^1f(x)
+3\nabla_x^2f(x)\otimes\nabla_x^2f(x) +(f(x)-t)\nabla_x^4f(x)
\end{array}
\eeqn
where $\otimes$ denotes the external product.

We finally establish the consequences of Theorem~\ref{general-th} in terms of the
functions involved in problem \req{genprob}.  Because this results makes
repeated used of Theorem 3.7 in \cite{CartGoulToin16a}, we first recall this
proposition.

\lthm{phi-res-th}{
\cite[Th.~3.7]{CartGoulToin16a} Suppose that $\psi$, a general objective function, is $q$ times continuously
differentiable and that $\nabla_x^q\psi$ is Lipschitz continous with constant
$L_{\psi,q}$ in an open neighbourhood of a point $x_\epsilon \in \calF$ of
radius larger than $\Delta_\epsilon$. Suppose also that, for some $\epsilon$,
\[
 \phi_{\psi,j}^{\Delta_\epsilon}(x_\epsilon) \leq \epsilon\Delta_\epsilon^j
 \tim{for} j=1,\ldots,q.
 \]
 Then
 \[
 \psi(x_\epsilon+d) \geq \psi(x_\epsilon) - 2 \epsilon \Delta^q
 \tim{ for all $d \in \calF(x_\epsilon)$ such that }
\|d\| \leq \left(\frac{q!\,\epsilon\Delta^q}{L_{\psi,q}}\right)^{\sfrac{1}{q+1}}.
 \]
}

\lbthm{lag-opt-th}{
  Suppose that AS.1--AS.4 hold and that, at
  $(x_\epsilon,t_\epsilon)$ and for some $\Delta_\epsilon > 0$,
  conditions \req{xinfeas} hold if $f(x_\epsilon)=t_\epsilon$ or conditions
  \req{mustat} hold for $\ii{q}$ if $f(x_\epsilon)>t_\epsilon$.
  \begin{enumerate}
  \item[(i) ] If $f(x_\epsilon) = t_\epsilon$ and, for  $j \in \ii{q}$,
    $\nabla_x^j \nu$ is Lipschitz continuous with constant $L_{\nu,j}$ in a
    neighbourhood of $x_\epsilon$ of radius larger than $\Delta_\epsilon$,
    then, for each $j \in \ii{q}$,
    \beqn{feas-lower}
    \|c(x_\epsilon)\| > \delta \ep
    \tim{ and }
    \|c(x_\epsilon+d)\| \geq \|c(x_\epsilon)\| - 2 \ed \|c(x_\epsilon)\| \Delta_k^j
    \eeqn
    for all $d \in \calF(x_\epsilon)$ such that
    \[
    \|d\| \leq \left(\bigfrac{j!\, \ed\|c(x_\epsilon)\|\Delta_\epsilon^j}
                       {L_{\nu,j}}\right)^{\sfrac{1}{j+1}}.
    \]
  \end{enumerate}
}{
  \begin{enumerate}
  \item[(ii) ] If $f(x_\epsilon) > t_\epsilon$, then, for
    \beqn{ydef}
    y_\epsilon = \frac{c(x_\epsilon)}{f(x_\epsilon)- t_\epsilon},
    \eeqn
    one has that
    \beqn{Lag-optim-j}
    \phi_{\Lambda,1}^{\Delta_\epsilon}(x_\epsilon,y_\epsilon)
    \leq \ed \Delta_\epsilon \|(1,y_\epsilon^T)\|
    \tim{and}
    \widehat{\phi}_{\Lambda,j}^{\Delta_\epsilon}(x_\epsilon,y_\epsilon)
    \leq \ed \Delta_\epsilon^j \|(1,y_\epsilon^T)\|
    \ms (j=2,3),
    \eeqn
    where $\widehat{\phi}_{\Lambda,j}^{\Delta_\epsilon}$ differs from
    $\phi_{\Lambda,j}^{\Delta_\epsilon}$ in that it uses the feasible 
    set $\calF(x_\epsilon) \cap \calM(x_\epsilon)$ instead of $\calF(x_\epsilon)$.
    Moreover, if $f$ and $c$ have Lipschitz continuous $j$-th derivatives with
    constants $L_{f,j}$ and $L_{c,j}$, respectively, then
    \beqn{obj-lower}
    \|c(x_\epsilon)\| \leq \delta \ep
    \tim{ and }
    f(x_\epsilon+d)
    \geq f(x_\epsilon) - 2 \ep \|y_\epsilon\| - 2 \ed \Delta_\epsilon^j \|(1,y_\epsilon^T)\|
    \eeqn
    for all $d$ such that $d\in\calM(x_\epsilon) \cap \calF(x_\epsilon)$
    whenever $j=2,3$, $\|c(x_\epsilon+d)\|\leq \epsilon$, and 
    \beqn{radj}
    \|d\| \leq \left(\bigfrac{j!\, \ed \Delta_\epsilon^j}
      {\sqrt{2}\max[L_{f,j},L_{c,j}]}\right)^{\sfrac{1}{j+1}},
    \eeqn
    Moreover, the second bound in \req{obj-lower} can be simplified to
    \beqn{simple-obj-lower}
    f(x_\epsilon+d)
    \geq f(x_\epsilon) - 2 \ed \Delta_\epsilon^j \|(1,y_\epsilon^T)\|
    \eeqn
    for any $d$ such that $d \in \calM(x_\epsilon) \cap \calF(x_\epsilon)$
    whenever $j=2,3$, \req{radj} holds, and for which  
    $c(x_\epsilon+d)=0$ or $c(x_\epsilon+d) = c(x_\epsilon)$.
  \end{enumerate}
}

\proof{
  Consider first the case where $f(x_\epsilon) = t_\epsilon$ (and thus
  $\|c(x_\epsilon)\| > \delta \ep$ because of
  Theorem~\ref{general-th}).  Note that we only need to consider the case
  where $\|c(x_\epsilon+d)\|\leq\|c(x_\epsilon)\|$. We have that, for 
  $d \in \calF(x_\epsilon)$,
  \[
  \|c(x_\epsilon+d)\|-\|c(x_\epsilon)\|
  = \frac{\|c(x_\epsilon+d)\|^2-\|c(x_\epsilon)\|^2}{\|c(x_\epsilon+d)\|+\|c(x_\epsilon)\|}
  \geq \frac{\nu(x_\epsilon+d)-\nu(x_\epsilon)}{\|c(x_\epsilon)\|}
  \]
  and the second part of \req{feas-lower} then follows from \req{xinfeas} and
  Theorem~\ref{phi-res-th} applied to the function $\nu$.

  Consider now the case where $f(x_\epsilon) > t_\epsilon$ (and thus
  $\|c(x_\epsilon)\| \leq \ep$ because of Theorem~\ref{general-th}).
  Focus first on the case where $j=1$. Theorem~\ref{general-th} then ensures that
  \[
  \phi_{\mu,1}^{\Delta_\epsilon}(x_\epsilon,t_\epsilon)
  \leq \ed \Delta_\epsilon \|r(x_\epsilon,t_\epsilon)\|.
  \]
  Using now \req{ydef} and
  \beqn{mu-L-1}
  \bigfrac{1}{f(x_\epsilon)-t_\epsilon}\nabla_x^1\mu(x_\epsilon,t_\epsilon)
   =  J(x_\epsilon)^T\bigfrac{c(x_\epsilon)}{f(x_\epsilon)-t_\epsilon}
        + \nabla_x^1f(x_\epsilon)
   =  J(x_\epsilon)^Ty_\epsilon+\nabla_x^1f(x_\epsilon)
   =  \nabla_x^1 \Lambda(x_\epsilon,t_\epsilon).
  \eeqn
  one has that \req{Lag-optim-j} holds for $j=1$.
  Moreover, applying Theorem~\ref{phi-res-th}, we obtain that
  \[
  \Lambda(x_\epsilon+d,y_\epsilon)
  \geq \Lambda(x_\epsilon,y_\epsilon)- 2 \ed \Delta_\epsilon \| (1,y_\epsilon^T)\|
  \]
  for all $d \in \calF(x_\epsilon)$ such that
  \[
  \|d\|
  \leq \sqrt{\frac{\|(1,y_\epsilon^T)\| \epsilon\Delta_\epsilon}{L_{f,1}+\|y_\epsilon\|L_{c,1}}}.
  \]
  Using now the fact that, for any $a \geq 0$, $\sqrt{2(1+a^2)}\geq 1+a$,
  we obtain that
  \beqn{yb1}
  \|(1,y_\epsilon^T)\| \geq \frac{1+\|y_\epsilon\|}{\sqrt{2}}.
  \eeqn
  Hence we deduce that, for all $d\in \calF(x_\epsilon)$ satisfying
  \beqn{rad1}
  \|d\|
  \leq \sqrt{\frac{\epsilon\Delta_\epsilon}{\sqrt{2}\max[L_{f,1},L_{c,1}]}},
  \eeqn
  we have that
  \beqn{f-lower-eps}
  f(x_\epsilon+d) + y_\epsilon^Tc(x_\epsilon+d)
  \geq f(x_\epsilon) + y_\epsilon^Tc(x_\epsilon)- 2 \ed \Delta_\epsilon \|  (1,y_\epsilon^T)\|
  \eeqn
  and hence, using the Cauchy-Schwarz inequality, that
  \[
  f(x_\epsilon+d)
  \geq f(x_\epsilon) - \|y_\epsilon\|\|c(x_\epsilon)-c(x_\epsilon+d)\|
       - 2 \ed \Delta_\epsilon \|(1,y_\epsilon)\|.
       \]
  If one additionally requests that $\|c(x_\epsilon+d)\| \leq \ep$, then, from
  the first part of \req{mustat},
  $
  \|c(x_\epsilon)-c(x_\epsilon+d)\| \leq 2 \ep
  $
  and therefore
  $
  f(x_\epsilon+d)
  \geq f(x_\epsilon) - 2 \ep \|y_\epsilon\| - 2 \ed \Delta_\epsilon \|(1,y_\epsilon^T)\|
  $
  for all $d \in \calF(x_\epsilon)$ such that \req{rad1} holds.
  Also note that, if $d$ exists such that $c(x_\epsilon+d)=0$, $x_\epsilon+d
  \in \calF$ and \req{rad1} holds, then \req{f-lower-eps} ensures that
  \beqn{fls1}
  f(x_\epsilon+d) \geq  f(x_\epsilon) - 2 \ed \Delta_\epsilon \|(1,y_\epsilon)^T\|
  \eeqn
  since $y_\epsilon^Tc(x_\epsilon) \geq 0$ because $f(x_\epsilon)- t_\epsilon > 0$.
  Similarly, if $d$ exists such that $c(x_\epsilon+d)=c(x_\epsilon)$, $d
  \in \calF(x_\epsilon)$ and \req{rad1} holds, then \req{f-lower-eps} ensures that
  \req{fls1} also holds.

  Now turn to the case where $f(x_\epsilon) > t_\epsilon$ and $j=2$.
  Observe now that, because of \req{D2M} and \req{Mdef}, 
  \beqn{mu-L-2}
  \nabla_x^2\Lambda(x_\epsilon,y_\epsilon)[d]^2
  = \frac{1}{f(x_\epsilon)- t_\epsilon}\nabla_x^2 \mu(x_\epsilon,t_\epsilon)[d]^2
  \tim{ for all }d \in \calM(x_\epsilon).
  \eeqn
  Now,  $\phi_{\mu,2}^{\Delta_\epsilon}(x_\epsilon) \leq \ed \Delta_\epsilon^2
  \|r(x_\epsilon,t_\epsilon)\|$  implies that
  \[
  \nabla_x^1 \mu(x_\epsilon,t_\epsilon)[d] + \half \nabla_x^2 \mu(x_\epsilon,t_\epsilon)[d]^2
  \geq - \ed \Delta_\epsilon^2 \|r(x_\epsilon,t_\epsilon)\|
  \]
  for all $d \in \calM(x_\epsilon) \cap \calF(x_\epsilon)$, and thus, dividing by
  $f(x_\epsilon)-t_\epsilon >0$ and using \req{mu-L-1} and \req{mu-L-2}, 
  \[
  \nabla_x^1 \Lambda(x_\epsilon,y_\epsilon)[d]
  + \half \nabla_x^2 \Lambda(x_\epsilon,y_\epsilon)[d]^2
  \geq - \ed \Delta_\epsilon^2 \|(1,y_\epsilon)\|
  \]
  for all $d \in \calM(x_\epsilon) \cap \calF(x_\epsilon)$.
  This in turn ensures that \req{Lag-optim-j} holds for $j=2$.
  Applying Theorem~\ref{phi-res-th} for the problem defining
  $\widehat{\phi}$, we deduce that
  \beqn{L-lower-e}
  \Lambda(x_\epsilon+d,y_\epsilon)
  \geq \Lambda(x_\epsilon,y_\epsilon)- 2 \ed \Delta_\epsilon^2 \|(1,y_\epsilon)\|
  \eeqn
  for all $d$ such that $d \in \calM(x_\epsilon) \cap \calF(x_\epsilon)$.  As
  a consequence, using \req{yb1} as above, we have that \req{obj-lower}
  holds for $j=2$ and all $d \in \calM(x_\epsilon)\cap \calF(x_\epsilon)$ such that
  \beqn{rad2}
  \|d\| \leq \left(\frac{2\ed\Delta_\epsilon^2}{\sqrt{2}\max[L_{f,2},L_{c,2}]}\right)^{\third}.
  \eeqn
  Applying the same reasoning as above, we deduce that
  \[
  f(x_\epsilon+d)
  \geq f(x_\epsilon) - 2 \ep \|y_\epsilon\| - 2 \ed \Delta_\epsilon^2 \|(1,y_\epsilon)\|
  \]
  if one additionally requests that $\|c(x_\epsilon+d)\| \leq \ep$.
  We may also, as for $j=1$, deduce from \req{L-lower-e} that
  $
  f(x_\epsilon+d)
  \geq f(x_\epsilon) - 2 \ed \Delta_\epsilon^2 \|(1,y_\epsilon)\|
  $
  for any $d$ such that $d \in \calM(x_\epsilon) \cap \calF(x_\epsilon)$
  and \req{rad2} holds and for which $c(x_\epsilon+d)=0$ or
  $c(x_\epsilon+d) = c(x_\epsilon)$.

  We finally turn to the case where $f(x_\epsilon) > t_\epsilon$ and $j=3$. It
  can be verified that, for $s_1 \in \calM(x_\epsilon)$,
  \beqn{mu-L-3}
  \begin{array}{lcl}
  \nabla_x^2\mu(x_\epsilon,t_\epsilon)[s_1,s_2]
  & = & \nabla_x^1c(x_\epsilon)[s_1].\nabla_x^1c(x_\epsilon)[s_2]+
     \nabla_x^1f(x_\epsilon)[s_1].\nabla_x^1f(x_\epsilon)[s_2]\\*[1ex]
  &   & + (f(x_\epsilon)-t_\epsilon) \nabla_x^2 \Lambda(x_\epsilon,y_\epsilon)[s_1,s_2]\\*[1ex]
  & = & (f(x_\epsilon)-t_\epsilon) \nabla_x^2 \Lambda(x_\epsilon,y_\epsilon)[s_1,s_2]
  \end{array}
  \eeqn
  and
  \beqn{mu-L-4}
  \begin{array}{lcl}
  \nabla_x^3\mu(x_\epsilon,t_\epsilon)[s_1]^3
   & = & 3 \Big[\sum_{i=1}^m \nabla_x^2c_i(x_\epsilon)[s_1]^2.\nabla_x^1c_i(x_\epsilon)[s_1]
     + \nabla_x^2f(x_\epsilon)[s_1]^2.\nabla_x^1f(x_\epsilon)[s_1]\Big]\\*[1ex]
   &   & + (f(x_\epsilon)-t_\epsilon) \nabla_x^3 \Lambda(x_\epsilon,y_\epsilon)[s_1]^3\\*[1ex]
   & = & (f(x_\epsilon)-t_\epsilon) \nabla_x^3 \Lambda(x_\epsilon,y_\epsilon)[s_1]^3.
   \end{array}
  \eeqn
  At termination we have that
  $\phi_{\mu,3}^{\Delta_\epsilon}(x_\epsilon) \leq \ed
  \Delta_\epsilon^3\|r(x_\epsilon,t_\epsilon)\|$,
  and thus, for all $d \in \calF(x_\epsilon)$,
  \[
  \nabla_x^1\mu(x_\epsilon,t_\epsilon)[d]+
  \half \nabla_x^2\mu(x_\epsilon,t_\epsilon)[d]^2 + 
  \sixth \nabla_x^3\mu(x_\epsilon,t_\epsilon)[d]^3
  \geq - \ed \Delta_\epsilon^3 \|r(x_\epsilon,t_\epsilon)\|.
  \]
  As for $j=1$ and $2$, and for every
  $d \in \calM(x_\epsilon)\cap \calF(x_\epsilon)$, the above relations imply
  that 
  \[
  \nabla_x^1\Lambda(x_\epsilon,y_\epsilon)[d]+
  \half \nabla_x^2\Lambda(x_\epsilon,y_\epsilon)[d]^2 + 
  \sixth \nabla_x^3\Lambda(x_\epsilon,y_\epsilon)[d]^3
  \geq - \ed \Delta_\epsilon^3 \|(1,y^T_\epsilon)\|,
  \]
  and therefore that \req{Lag-optim-j} holds for $j=3$.
  Applying Theorem~\ref{phi-res-th} again, we now deduce that
  \[
  \Lambda(x_\epsilon+d,y_\epsilon)
  \geq \Lambda(x_\epsilon,y_\epsilon)- 2 \ed \Delta_\epsilon^2 \|(1,y_\epsilon)\|
  \]
  for all $d \in \calM(x_\epsilon) \cap \calF(x_\epsilon)$.
  As for the previous cases, this implies that \req{obj-lower} holds for for
  $j=3$ and, using \req{yb1} once more, for all $d \in \calM(x_\epsilon)\cap
  \calF(x_\epsilon)$ satisfying
   \beqn{rad3}
  \|d\| \leq \left(\frac{6\ed\Delta_\epsilon^3}{\sqrt{2}\max[L_{f,3},L_{c,3}]}\right)^{\third}.
  \eeqn
  The inequality \req{simple-obj-lower} is obtained as for the cases where $j=1,2$.
} 

\noindent
We verify that \req{Lag-optim-j} for $j=1$
is the scaled first-order criticality condition considered in
\cite{CartGoulToin13a} (Theorem~\ref{lag-opt-th} thus subsumes the
analysis presented in that reference) and is equivalent to
\[
\| P_{\calT_(x_\epsilon)}[-\nabla_x^1\Lambda(x_\epsilon,y_\epsilon)] \| \leq \ed \Delta_\epsilon  \| (1,y_\epsilon^T)\|,
\]
which corresponds to a scaled version of the first-order criticality condition
considered in \cite{BirgGardMartSantToin16}.

\subsection{Beyond third-order optimality?}

We have now proved that, if an approximate $q$-th order critical points for the
convexly constrained problem can be obtained by an inner algorithm
at a given evaluation complexity, then the same result holds
for the critical points of $\|c(x)\|$ whenever \outeralg\ 
terminates with an infeasible stationary point of the constraint violation
(either at Phase~1 or at \req{xinfeas}).  When \outeralg\ 
terminates with \req{mustat}, we have shown in Theorem \ref{lag-opt-th} that
similar results hold for criticality of orders one, two and three.

As indicated already, the situation becomes considerably more complicated for
higher orders.  The first difficulty, which we covered in Section~\ref{cnc-s},
is that the conditions \req{mgradLinN}-\req{ccnc2-e} involve, for higher
orders, the geometry of the feasible arcs in a way which is hard to exploit.
Moreover, the fact that we could
derive, in Theorem~\ref{lag-opt-th}, some lower bounds on the objective
function values by exploiting information at orders one up to three is strongly
dependent of the observation that, in the suitable subspace,
\beqn{mu-L}
\frac{1}{f(x_\epsilon)-t_\epsilon }\nabla_x^j\mu(x_\epsilon,t_\epsilon)
= \nabla_x^j\Lambda(x_\epsilon,y_\epsilon)
\ms (j = 1,2,3)
\eeqn
(see \req{mu-L-1}, \req{mu-L-2}, \req{mu-L-3} and \req{mu-L-4}), which in turn
ensures that minimizing $\mu(x,t)$ with respect to $x$ on the said subspace also results in
minimizing $\Lambda(x,y)$ with respect to $x$ on the same subspace\footnote{For order
three, it is fortunate that terms in \req{mu-L-3} and \req{mu-L-4} involving the
second derivatives always appear in product with terms involving the first,
which is the reason why the minimization subspace at order three is not smaller
than that at order two.}. Is this crucial property maintained for high orders?
We now 
show that the answer to this question is negative
for orders four and beyond, due to the ever more distant relationship between
$\nabla_x^j\mu(x_\epsilon,t_\epsilon)$ and $\nabla_x^j \Lambda(x_\epsilon,y_\epsilon)$
when $j$ grows, which is apparent when considering the expressions
\req{D1M}-\req{D4M}.  Indeed, the terms   
\beqn{problematic-terms}
\begin{array}{lcl}
\lefteqn{\bigfrac{3}{f(x_\epsilon)-t_\epsilon}
\left[\bigsum_{i=1}^m (\nabla_x^2c_i(x)\otimes\nabla_x^2c_i(x))[d]^4
  + (\nabla_x^2f(x)\otimes\nabla_x^2f(x))[d]^4\right]} \\
\ms\ms\ms\ms\ms\ms\ms\ms
& = & \bigfrac{3}{f(x_\epsilon)-t_\epsilon}
\left[\bigsum_{i=1}^m \Big(\nabla_x^2c_i(x)[d]^2\Big)^2
  + \Big(\nabla_x^2f(x)[d]^2 \Big)^2\right]
\end{array}
\eeqn
in \req{D4M} would only vanish in general if $d \in \ker^2[\nabla_x^2f(x)]
\cap \ker^2[\nabla_x^2c(x)]$.  Although this is formally reminiscent of the definition
of $\calM(x)$ in \req{Mdef}, this crucial inclusion now no longer follows from
lower-order conditions.

\begin{figure}[htbp]
\vspace*{1.5mm}
\includegraphics[height=4.3cm]{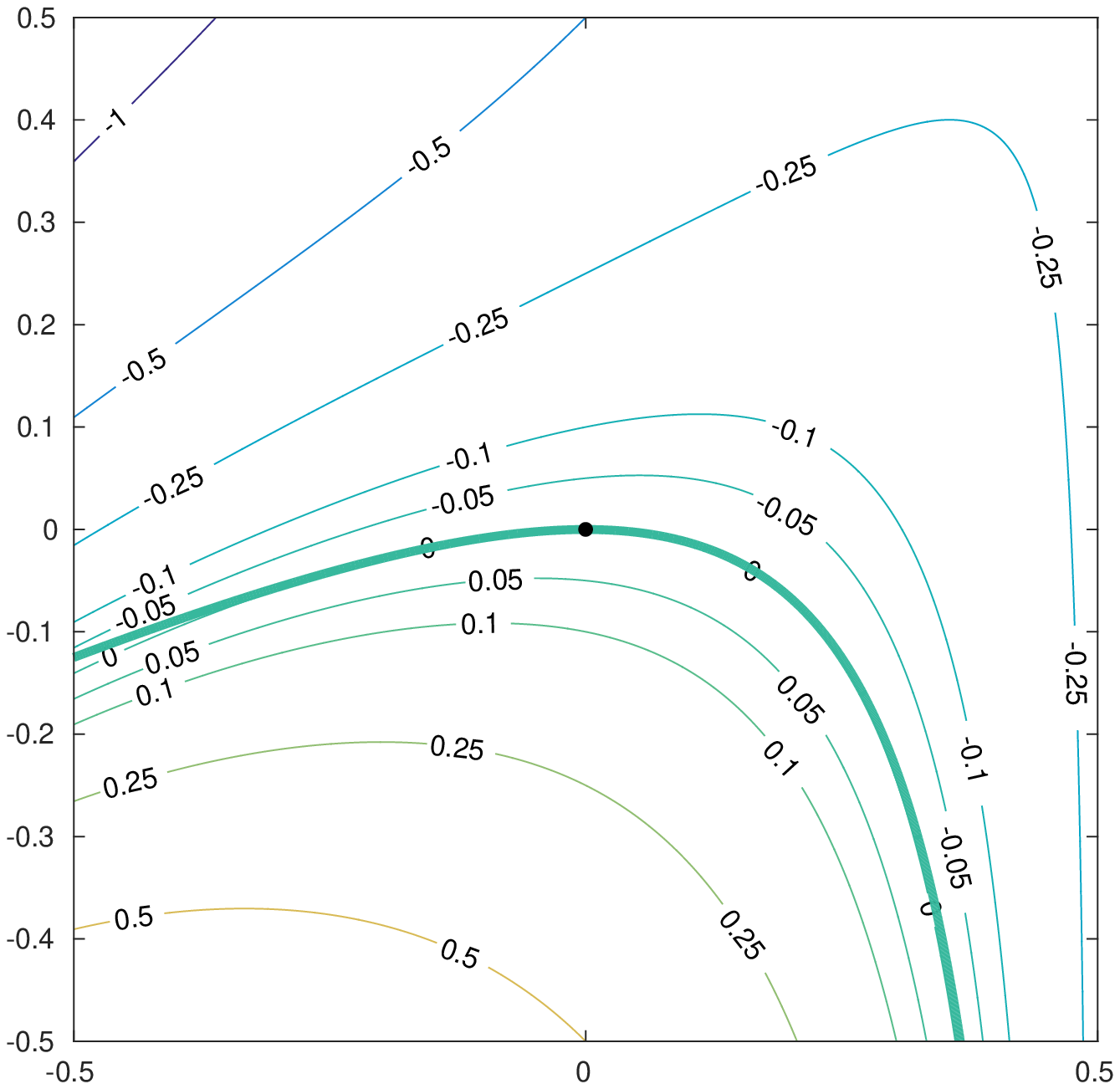}
\hspace*{9mm}
\includegraphics[height=4.3cm]{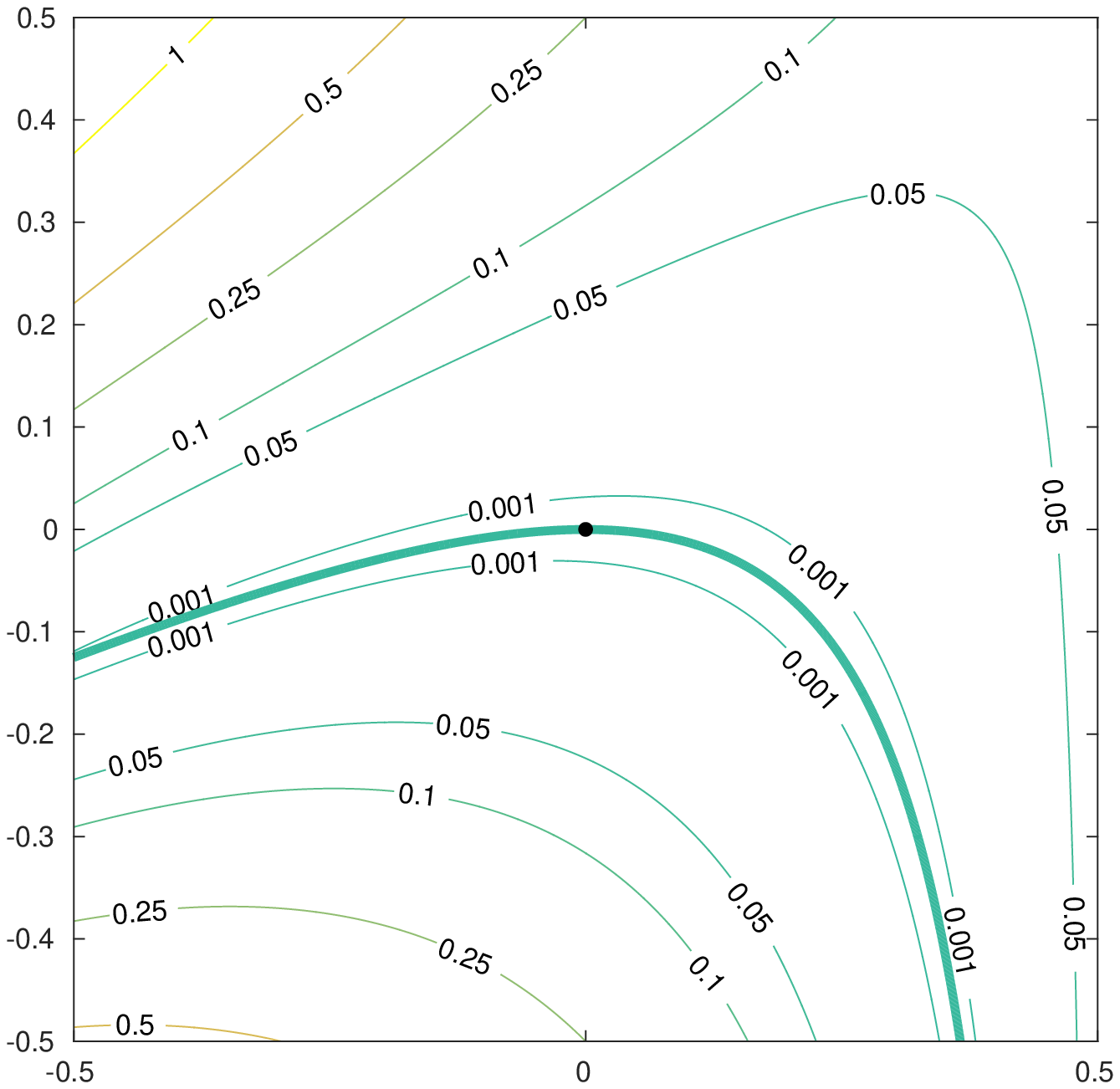}  
\hspace*{9mm}
\includegraphics[height=4.3cm]{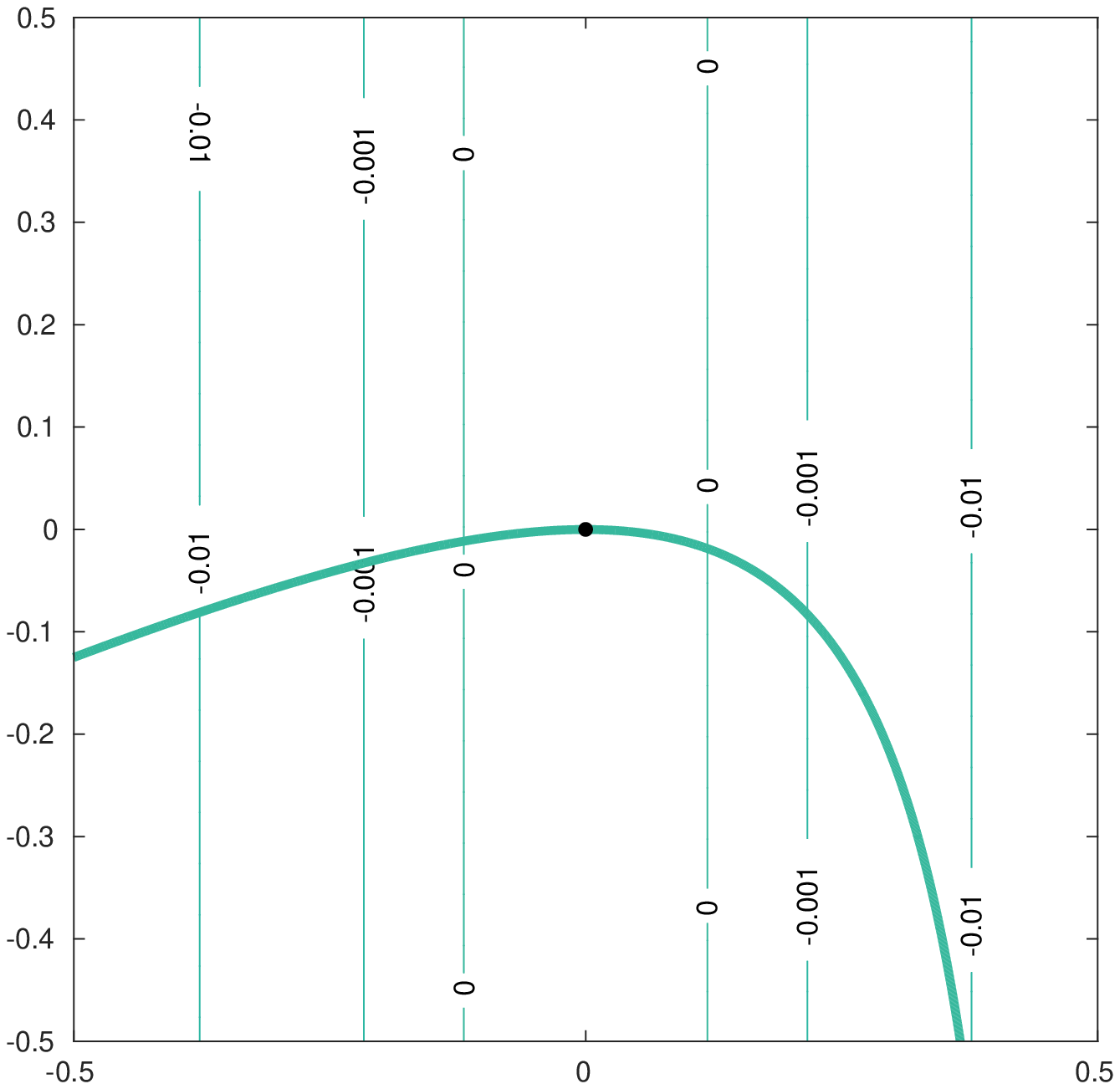} 
\caption{\label{theprob-fc}Contour lines for \req{theprob} (left), $\mu(x,t_\epsilon)$
  (center) and $\Lambda(x,y_\epsilon)$ (right)with the
  constraint shown as a thick curve} 
\end{figure}

This is illustrated by what happens on the problem
\beqn{theprob}
\min_{x_1,x_1} - x_2 - x_1^2 + x_1x_2 - \half x_1^4
\tim{ subject to }
\varepsilon + x_2 + x_1^2 - x_1 x_2 = 0
\eeqn
for some $\varepsilon \in (0,1]$.
If we consider $x_\epsilon = (0,0)$ and $t_\epsilon = -\varepsilon$
(yielding $y_\epsilon= 1$), then one can verify (see Appendix)
that $\mu(0,t_\epsilon)$ satisfies the necessary conditions for a fourth order
minimizer at the origin while the problem itself has a global (fourth order)
constrained maximizer. Figure~\ref{theprob-fc}
shows the contour lines of the objective function with the constraint set
superimposed as a thick curve (left),  the contour lines of
$\mu(x,t_\epsilon)$ (center) and $\Lambda(x,y_\epsilon)$ (right).

It is worthwhile to note that the above discussion has wider
implications. Indeed the first of the problematic terms in
\req{problematic-terms} not only occurs in the function $\mu(x,t)$ used in
this paper, but also when applying to problem \req{genprob} a quadratic,
$\ell_1$ or $\ell_\infty$ penalty function, a classical augmented
Lagrangian approach, or a sequential quadratic programming method using a
merit function depending on such penalty terms. The same difficulty may also
occur if more general penalizations of the type $p(\nu(x))$ (for some
increasing smooth function $p$ from $\Re^+$ to $\Re^+$) are employed. Indeed,
consider the derivatives of $p(\nu(x))$.  One verifies that
\[
\begin{array}{lcl}
\nabla_x^4 p(\nu(x))
&=& p''''(\nu(x)) [\nabla_x^1\nu(x)]^{4 \otimes}
     + 6 p'''(\nu(x))
    \nabla_x^2\nu(x) \otimes [\nabla_x^1\nu(x)]^{2 \otimes}\\*[2ex]
& & + 4 p''(\nu(x))
      \nabla_x^3\nu(x) \otimes \nabla_x^1\nu(x)
     + 3 p''(\nu(x))
      [\nabla_x^2\nu(x)]^{2\otimes}\\*[2ex]
& & + p'(\nu(x)) \nabla_x^4 \nu(x)\\*[3ex]
\end{array}
\]
whose last term, together with
\[
\nabla_x^4 \nu(x)
= \sum_{i=1}^m \Bigg[4\,\nabla_x^3c_i(x)\otimes\nabla_x^1c_i(x)
  +3\,[\nabla_x^2c_i(x)]^{2\otimes} +c_i(x)\nabla_x^4c_i(x)\Bigg],
\]
indicates that the troublesome terms involving $[\nabla_x^2c_i(x)]^{2\otimes}$
do not vanish unless $p'(\nu(x))$ also vanishes with $\nu(x)$.

\emph{None of the linear or quadratic penalization approaches can therefore be
  expected to reliably produce critical points of orders four or more.}
Innovative techniques are thus needed if one is interested to compute
high-order critical points of \req{genprob} of higher order.  One
possible research direction is to follow the propositions formulated in
\cite{CartGoulToin15c} and to exploit penalization terms of order higher than
two in the definitions of $\nu$ and $\mu$, for which an improved evaluation
complexity bound is already available for the subproblem solution.

\numsection{Conclusions and discussion}\label{disc-s}

We have formulated and analyzed, in Section~\ref{cnc-s}, the necessary
conditions for high-order optimality in nonlinear optimization problems
involving both convex set constraints and nonlinear equalities. We have also
discussed the difficulties inherent to their form for third-order critical
points and higher.

We then have shown in Sections~\ref{general-s} and \ref{complexity-s} that the
evaluation complexity of finding an approximate $q$-th-order scaled critical
point ($q=1,2,3$) for a large class of smooth nonlinear optimization problem
involving both equality and inequality constraints is at most
$O(\ep^{1-\pi}\ed^{-\pi})$ evaluations of the objective function, constraints
and their derivatives, where $\epsilon^\pi$ is the order of the guaranteed
objective function decrease during the successful iterations of an underlying
inner algorithm for convexly constrained least-squares problems. We refer here
to an ``approximate scaled critical point'' in that such a point is required
to satisfy \req{xinfeas} or \req{mustat}, where the accuracy is scaled by the
size of the constraint violation or that of the Lagrange multipliers.  In
particular, the above results provide the first evaluation complexity bound
for second- and third-order criticality in the case involving general
inequality and equality constraints.

This result also corrects an unfortunate error\footnote{The second equality in
the first equation of Lemma~3.4 in \cite{CartGoulToin13a} only holds if one
is ready to flip the gradient's sign if necessary.} in the first-order
analysis of \cite{CartGoulToin13a}, that allows a vector of Lagrange
multipliers whose sign is arbitrary (in line with a purely first-order setting
where minimization and maximization are not distinguished).  The present
analysis now yields the multiplier with the sign associated with minimization.

Interestingly, an $O(\ep \ed^{-(p+1)/p} \min[\ed,\ep]^{-(p+1)/p})$ evaluation
complexity bound was also proved by Birgin, Gardenghi, Mart\'{i}nez, Santos
and Toint in \cite{BirgGardMartSantToin16} for first-order \emph{unscaled},
standard  KKT  conditions and in the least expensive of three cases depending
on the degree of degeneracy identifiable by the algorithm\footnote{This result
also assumes boundedness of $f(x_1)$.}. Even if the bounds for the scaled and
unscaled cases coincide in order when $\ep \leq \ed$, comparing the two
results for first-order critical points is not straightforward.  On
one hand the scaled conditions take into account the possibly different
scaling of the objective function and constraints.  On the other hand the same
scaled conditions may result in earlier termination with \req{mustat} if the
Lagrange multipliers are very large, as \req{mustat} is then consistent with the
weaker requirement of finding a John's point. But the framework discussed in
the present paper also differs from that of \cite{BirgGardMartSantToin16} in
additional significant ways.  The first is that second-order critical points
are now covered in the analysis.  If we now restrict the scope to first-order,
the present paper provides a potentially stronger version of the termination
of the algorithm at infeasible points (in Phase~1): indeed the second part of
\req{xinfeas} can be interpreted as requiring that the size of the feasible
gradient of $\|c(x)\|$ is below $\ed$, while \cite{BirgGardMartSantToin16}
considers the gradient of $\|c(x)\|^2$ instead. The second is that, if
termination occurs in Phase~2 for an $x_\epsilon$ such that
$\phi_{\nu,1}^{\Delta_k}(x_\epsilon)$ is of order $\ep\ed\Delta_k$ (thereby
covering the case where $f(x_\epsilon)=t_k$ discussed in
Theorem~\ref{general-th}) and $x_\epsilon \in \calF^0$, then 
$\|P_{\calT_*}[-\nabla_x^1\nu(x_*)]\|=\|\nabla_x^1\nu(x_*)\|$ is of
the same order and Birgin \textit{et al.} show that, in
this case, the {\L}ojaciewicz inequality \cite{Loja65} must fail for $c$ in
the limit for $\ep$ and $\ed$ tending to zero (see
\cite{BirgGardMartSantToin16} for details).  This observation is interesting
because smooth functions satisfy the {\L}ojaciewicz inequality under
relatively weak conditions, implying that termination in these circumstances
is unlikely. The same information is also obtained in
\cite{BirgGardMartSantToin16}, albeit at the price of worsening the evaluation
complexity bound mentioned above by an order of magnitude in $\ed$. We also
note that the approach of \cite{BirgGardMartSantToin16} requires the
minimization, at each iteration, of a residual whose second derivatives are
discontinuous, while all functions used in the present paper are $p$ times
continuously differentiable.  A final difference between the two approaches is
obviously our introduction of $\phi_{\mu,j}^{\Delta_k}$ in the expression of
the criticality condition in Theorem~\ref{general-th} for taking the
inequality constraints into account.

Will regularization-based methods provide better evaluation complexity
bounds when using polynomial models of higher order? Can the limitations of
penalty approaches for finding high-order solutions for equality constrained
problems be circumvented? These and many other questions remain open at this stage.

\section*{\footnotesize Acknowledgements}

{\footnotesize

The authors would like to thank Oliver Stein for suggesting reference
\cite{Hoga73}. The work of the second author was supported by EPSRC grant
EP/M025179/1. The third author acknowleges the support
provided by the Belgian Fund for Scientific Research (FNRS), the Leverhulme
Trust (UK), Balliol College (Oxford, UK), the Department of Applied
Mathematics of the Hong Kong Polytechnic University, ENSEEIHT (Toulouse,
France) and INDAM (Florence, Italy).



}

\vspace*{4mm}
\appendix
\noindent
\textbf{\Large Appendix}

\appnumsection{Details of the derivation of \req{fx1-bound}}\label{A-fx1-bound}

\noindent
For the trust-region algorithm,
\[
f(x_1^*)
\leq  f(x_1) +
\left(\bigmax_{\xi \in \cup_{j\in \calS}[x_j,x_{j+1}]}\|\nabla_x^1\nu(\xi)\| \right)
    \Delta_{\max} \left(\left\lfloor \kap{CC}^{\|c\|} \|c(x_1)\| \,
    \ep^{-q}\ed^{-(q+1)} \right\rfloor
    + 1 \right).
\]
For the regularization algorithm,
\[
\begin{array}{lcl}
f(x_1)
& \leq  & f(x_0) + \left(\bigfrac{\nu(x_0)(p+1)!}{\eta
  \sigma_{\min}}\right)^{\frac{1}{p+1}} \times \\*[1ex]
& & \bigmax_{\xi \in \cup_{j\in \calS}[x_j,x_{j+1}]}\|\nabla_x^1\nu(\xi)\|
    \left\{\left\lfloor \kap{CC}^{\|c\|} \|c(x_0)\| \,
    \bigmax\left[\ep^{-1}, \ep^{-\frac{1}{p}}\ed^{-\frac{p+1}{p+1-q}}\right] \right\rfloor
    + 1 \right\}.
\end{array}
\]

\appnumsection{Details for the example \req{theprob}}\label{A-theprob}

We prove the validity of the statement made after the definition of problem
\req{theprob}, namely that $\mu(0,t_\epsilon)$
satisfies the necessary conditions for a fourth order minimizer at the origin while
the problem itself has a global (fourth order) constrained maximizer.

Let $T_3(x) = x_2 + x_1^2 - 2x_1x_2$ and define, for some $\varepsilon \in (0,1]$,
\beqn{mu4l-fc}
f(x) =  - T_3(x) - \half x_1^4
\tim{ and }
c(x) =  \varepsilon + T_3(x).
\eeqn
and thus, for a given multiplier $y$, 
\beqn{mu4l-L}
\Lambda(x,y)
= -T_3(x) - \half x_1^4  + y[ \varepsilon + T_3(x) ]
\eeqn
We have that
\beqn{mu4l-derT}
\nabla_x^1 T_3(x) = \cvect{ 2x_1 -2x_2\\ 1 -2x_1 },
\ms
\nabla_x^2 T_3(x) = \mat{rr}{ 2 & -2 \\ -2 & 0 }
\tim{and}
\nabla_x^3T_3(x) = 0.
\eeqn
Thus, at the origin and for $t_\epsilon = - \varepsilon$
\beqn{cTf}
c(0) = \varepsilon = f(0) - t_\epsilon
\tim{ and }
\nabla_x^jc(0) = \nabla_x^jT_3(0) = - \nabla_x^jf(0)
\tim{ for }
j = 1, 2, 3.
\eeqn
As a consequence, the choice $y=1$, \req{mu4l-L} and \req{mu4l-derT} ensure that
$
\Lambda(x,1)
= \varepsilon  - \half x_1^4
$
as well as
\beqn{mu4l-derL0}
\nabla_x^1 \Lambda(0,1) = 0,
\ms
\nabla_x^2 \Lambda(0,1) = 0,
\ms
\nabla_x^3 \Lambda(0,1) = 0,
\ms
\nabla_x^4 \Lambda(0,1) = \nabla_x^4f(0)=-12 e_1^{\otimes 4}.
\eeqn
Using \req{D1M}-\req{D4M} and \req{cTf}, we also have that, for $t=-\varepsilon$, 
\beqn{d1mu}
\nabla_x^1\mu(0,t_\epsilon)
= (c(0) - f(0) + t_\epsilon) \nabla_x^1 T_3(0)
= ( \varepsilon - 0 -\varepsilon ) e_2
= 0,
\eeqn
\beqn{d23mu}
\nabla_x^2\mu(0,t_\epsilon)
= 2 \, \nabla_x^1T_3(0)\otimes \nabla_x^1T_3(0)
= 2 e_2e_2^T,
\ms
\nabla_x^3 \mu(0,t_\epsilon)
= 6 \, \nabla_x^2T_3(0)\otimes\nabla_x^1T_3(0)
= 0
\eeqn
and, using the last equation in \req{mu4l-derL0},
\beqn{d4mu}
\begin{array}{lcl}
\nabla_x^4 \mu(0,t_\epsilon)
& = & 6 \,\nabla_x^2T_3(0)\otimes\nabla_x^2T_3(0)
+ c(0)\nabla_x^4c(0) + (f(0)-t_\epsilon)\nabla_x^4f(0) \\*[1ex]
& = & 12 \left[ \mat{rr}{ 1 & -1 \\ -1 & 0 }^{\otimes 2}
      - \varepsilon e_1^{\otimes 4}\right].
\end{array}
\eeqn
(Notice the contribution of the first term in the bracketed
expression, potentially dwarfing that of the second for sufficiently small
$\epsilon$.)  

Let us attempt to verify \req{mgradLinN}-\req{ccncf} with $q=4$ for the
problem of minimizing $\mu(x,t_\epsilon)$ with $s_1 \in \ker^2[\nabla^2\mu] =
\spanset{e_1}$. We have that \req{mgradLinN} holds because of \req{d1mu}. We
also obtain, from \req{d1mu}-\req{d4mu}, that, for $s_1 = \tau e_1$  for some
$\tau \in \Re$ and for any choice of $s_2, s_3, s_4 \in \Re^n$,
\[
\nabla_x^2\mu(0,t_\epsilon)[s_2] + \half \nabla_x^2\mu(0,t_\epsilon)[s_1]^2
= 0^Ts_2 + \tau e_2 e_2^Te_1
= 0,
\]
\[
\nabla_x^3\mu(0,t_\epsilon)[s_3] + \tau \nabla_x^2\mu(0,t_\epsilon)[s_1,s_2]
+ \frac{\tau^3}{6} \nabla_x^3\mu(0,t_\epsilon)[s_1]^3
= 0^Ts_3 + \tau s_2^Te_2 e_2^T e_1 + \frac{\tau^3}{6} 0[s_1]^3
= 0
\]
and
\beqn{mu4l-cond4mu}
\begin{array}{ll}
\nabla_x^1&\!\!\!\!\!\mu(0,t_\epsilon)[s_4] + \nabla_x^2\mu(0,t_\epsilon)[s_1,s_3]
+ \half\nabla_x^2\mu(0,t_\epsilon)[s_2]^2+\half\nabla_x^3\mu(0,t_\epsilon)[s_1,s_1,s_2]
+ \sfrac{1}{24}\nabla_x^4\mu(0,t_\epsilon)[s_1]^4 \\*[1.5ex]
&= 0^Ts_4 + + \tau s_3^Te_2 e_2^Te_1 + \half (e_2^Ts_2)^2
    + \frac{\tau^2}{2} 0^T[e_1,e_1,s_2]\\*[1.5ex]
& \ms\ms\ms +\sfrac{12}{24}\left[\left\|e_1^T\mat{rr}{1 & -1\\ -1 & 0}e_1\right\|^2
  -\varepsilon \right] \tau^4 (e_1^Te_1)^4 \\*[2.5ex]
&= \half (e_2^Ts_2)^2+ \sfrac{1}{2}(1-\varepsilon ) \tau^4(e_1^Te_1)^4.
\end{array}
\eeqn
The choice of $s_2$, $s_3$ and $s_4$ is however constrained by \req{ccncf} for
$i=1,2,3$, in that those vector must also satisfy the equations
\[
\nabla_x^1c(0)[s_1] = 0 = \tau e_2^Te_1,
\]
\[
\nabla_x^1c(0)[s_2] +\half \nabla_x^2c(0)[s_1]^2= 0 = e_2^Ts_2
+ \tau^2 e_1^T\mat{rr}{ 1 & -1 \\ -1 & 0 }e_1,
\]
\[
\nabla_x^1c(0)[s_3] + \nabla_x^2c(0)[s_1,s_2] + \half \nabla_x^3c(0)[s_1]^3
= 0 = e_2^Ts_3 + 2\tau(e_1-e_2)^Ts_2 + \tau^3 0^T[e_1]^3.
\]
and
\[
\begin{array}{ll}
\nabla_x^1&\!\!\!\!\!c(0)[s_4] + \nabla_x^2c(0)[s_1,s_3] + +\half\nabla_x^2c(0)[s_2]^2
+ \half \nabla_x^3c(0)[s_1,s_1,s_2]^2 + \sfrac{1}{24}\nabla_x^4c(0)[s_1]^4\\*[1ex]
&= 0 = e_2^Ts_4 + e_1^Ts_2(e_1^Ts_2-2\tau^2) + 2 \tau e_1^Ts_3 - 4\tau^2.
\end{array}
\]
The second, third and fourth of these conditions impose constraints on the values of
$e_2^Ts_2$, $e_2^Ts_3$ and $e_2^Ts_4$. In particular, the second implies that
$e_2^Ts_2= -\tau^2$, which we may then substitute in \req{mu4l-cond4mu} and
deduce that
\beqn{n4mu}
\begin{array}{ll}
\nabla_x^4&\!\!\!\!\!\mu(0,t_\epsilon)[s_4] + \nabla_x^2\mu(0,t_\epsilon)[s_1,s_3]
+ \half\nabla_x^2\mu(0,t_\epsilon)[s_2]^2+\half\nabla_x^3\mu(0,t_\epsilon)[s_1,s_1,s_2]
+\sfrac{1}{24}\nabla_x^4\mu(0,t_\epsilon)[s_1]^4 \\*[1.5ex]
& = \half \tau^4 + (\half - \varepsilon) \tau^4
  = ( 1 - \half \varepsilon) \tau^4
  \geq 0.
\end{array}
\eeqn
We therefore obtain that, for all $\varepsilon \in (0,1]$, $x_*$ satisfies
the necessary conditions of Theorem~\ref{ccnc-th} with $q=4$, except that
$c(x_*) = \varepsilon$. However \req{mu4l-derL0} shows that
$\Lambda(x,y_\epsilon)$ is a polynomial of degree 4 with a global maximizer at
the origin, independently of the value of $\varepsilon$. Letting $\varepsilon$
tend to zero and using the fact that all quantities in the example depend
continuously on this parameter then allows to conclude.  

\end{document}